\documentclass[11pt,a4j]{article}
\usepackage{amsmath,amsthm,amssymb,bm}
\usepackage[authoryear,round]{natbib}
\usepackage{graphicx}
\usepackage{algorithm}
\usepackage{algorithmic}
\usepackage{color, xcolor}
\usepackage{breqn}
\usepackage{ulem}

\usepackage{indentfirst}
\usepackage{url}
\def\X{{\bm{X}}}

\def\F{{\bm{F}}}
\def\eps{{\bm{\varepsilon}}}
\def\0{{\bm{0}}}

\def\bfmu{\bm \mu}

\definecolor{darkgreen}{rgb}{0.0, 0.5, 0.0}

\newtheoremstyle{italicTheorem}% name
  {3pt}% Space above
  {3pt}% Space below
  {\itshape}% Body font
  {}% Indent amount
  {\itshape}% Theorem head font
  {.}% Punctuation after theorem head
  {.5em}% Space after theorem head
  {}% Theorem head spec (can be left empty, meaning ‘normal’)

\theoremstyle{plain}
\newtheorem{thm}{Theorem}
\newtheorem{pro}{Proposition}
\newtheorem*{cor*}{Corollary}
\newtheorem{cor}{Corollary}
\newtheorem{lem}{Lemma}
\newtheorem*{thm*}{Theorem}

\theoremstyle{italicTheorem}
\newtheorem*{rem*}{Remark}
\newtheorem{rem}{Remark}
\theoremstyle{definition}
\newtheorem{dfn}{Definition}
\newtheorem*{dfn*}{Definition}

\allowdisplaybreaks[1]

\begin{document}

\begin{center}
    \textbf{\Large Algebraic Approach for Orthomax Rotations}
\end{center}
\begin{center}
\large {Ryoya Fukasaku$^{1}$, Michio Yamamoto$^{2,5,6}$, Yutaro Kabata$^{3}$, Yasuhiko Ikematsu$^{4}$, Kei Hirose$^{4}$}
\end{center}

\begin{flushleft}
{\footnotesize
$^1$ Faculty of Mathematics, Kyushu University, 744 Motooka, Nishi-ku, Fukuoka 819-0395, Japan \\
$^2$ Graduate School of Human Sciences, Osaka University, 1-2 Yamadaoka, Suita, Osaka 565-0871, Japan\\
$^3$ School of Information and Data Sciences, Nagasaki University, Bunkyocho 1-14, Nagasaki 852-8131, Japan \\
$^4$ Institute of Mathematics for Industry, Kyushu University, 744 Motooka, Nishi-ku, Fukuoka 819-0395, Japan \\
$^5$ RIKEN Center for Advanced Intelligence Project, 1-4-1, Nihonbashi, Chuo-ku, Tokyo 103-0027, Japan \\
$^6$ Data Science and AI Innovation Research Promotion Center, Shiga University, 1-1-1, Banba, Hikone, Shiga, 522-8522, Japan \\
\vspace{1.2mm}}
{\it {\small E-mail: fukasaku@math.kyushu-u.ac.jp, yamamoto.michio.hus@osaka-u.ac.jp, kabata@nagasaki-u.ac.jp, ikematsu@imi.kyushu-u.ac.jp, hirose@imi.kyushu-u.ac.jp}}	
\end{flushleft}

\vspace{1.5mm}

\begin{abstract}
  In exploratory factor analysis, rotation techniques are employed to derive interpretable factor loading matrices. Factor rotations deal with equality-constrained optimization problems aimed at determining a loading matrix based on measure of simplicity, such as ``perfect simple structure'' and ``Thurstone simple structure.'' Numerous criteria have been proposed, since the concept of simple structure is fundamentally ambiguous and involves multiple distinct aspects. However, most rotation criteria may fail to consistently yield a simple structure that is optimal for analytical purposes, primarily due to two challenges. First, existing optimization techniques, including the gradient projection descent method, exhibit strong dependence on initial values and frequently become trapped in suboptimal local optima. Second, multifaceted nature of simple structure complicates the ability of any single criterion to ensure interpretability across all aspects. In certain cases, even when a global optimum is achieved, other rotations may exhibit simpler structures in specific aspects. To address these issues, obtaining all equality-constrained stationary points — including both global and local optima — is advantageous. Fortunately, many rotation criteria are expressed as algebraic functions, and the constraints in the optimization problems in factor rotations are formulated as algebraic equations. Therefore, we can employ computational algebra techniques that utilize operations within polynomial rings to derive exact all equality-constrained stationary points. Unlike existing optimization methods, the computational algebraic approach can determine global optima and all stationary points, independent of initial values. We conduct Monte Carlo simulations to examine the properties of the orthomax rotation criteria, which generalizes various orthogonal rotation methods. 
\end{abstract}

\section{Introduction}

The factor analysis model is a latent variable model initially introduced by \citep{Spearman.1904}. Recently, its applications have expanded social and behavioral sciences to encompass diverse fields, including marketing, life sciences, materials sciences, and energy sciences \citep{Lin.2019,Shkeer.2019,Shurrab.2019,Kartechina.2020,Vilkaite-Vaitone.2022}. A key objective in factor analysis is the estimation of the loading matrix, which represents the influence of latent variables — termed common factors — on observed variables.

This matrix is typically estimated via maximum likelihood or least squares method. After that, in order to enhance interpretability, factor rotations aim to achieve a ``simple structure," quantified through various measures. Notably, ``perfect simple structure'' and ``Thurstone simple structure'' are widely known as such a simple structure. 

Factor rotations deal with equality-constrained optimization problems, aiming to maximize or minimize specific rotation criterion. The concept of simple structure is inherently ambiguous, encompassing multiple distinct aspects. Thus, numerous rotation criteria have been proposed. However, most existing rotation criteria may fail to consistently yield a simple structure conducive to interpretability due to two primary challenges. First, optimization methods such as the gradient projection descent method exhibit strong dependence on initial values owing to the nonlinearity of rotation criteria. Therefore, such methods may result in suboptimal local optima.  Many software packages employ a single initial value, which may obscure suboptimal solutions from analysts unfamiliar with this dependence.

Second, given the multifaceted nature of simple structure, no single criterion can ensure that a loading matrix is interpretable across all aspects. Analysts select criteria based on specific simplicity objectives or empirical considerations. However, in certain cases, other rotations, e.g. stationary points, may yield a simpler structure compared to the global optimum in specific aspects. As conventional methods typically yield only a single solution, analysts are constrained from exploring potentially more informative alternatives.

This limitation is intrinsic because conventional approaches predominantly yield a single solution. To address these constraints, obtaining all equality-constrained stationary points — including both global and local optima — would enable a more informed selection of the most interpretable factor loadings. Fortunately, most of criteria can be formulated as algebraic functions, and the constraints in the optimization problems in factor rotations are formulated as algebraic equations. Since algebraic equations can be derived by applying the method of Lagrange multipliers to an equality-constrained optimization problem, the problem can be solved by solving these algebraic equations. Even more fortunately, algebraic operations can yield all solutions to such algebraic equations.

For instance, \citep{Jennrich2004} conducted simulations on orthogonal rotation criteria using the gradient projection algorithm for orthomax rotations introduced by \citep{Jennrich2001}. However, this method is highly sensitive to initial values and is incapable of computing all possible solutions. To address this issue, the present study eliminates such dependency by employing simulations based on algebraic operations. Since existing research has not utilized this approach, as will be demonstrated later, it yields novel findings in this study.

This study examines the ``orthomax criteria'' as introduced in \citep[section 4.6]{Harman}, which generalizes various orthogonal rotation criteria, including the ``quartimax criterion,'' ``varimax criterion,'' ``equamax criterion,'' ``parsimax criterion'', and etc. Specifically, orthomax rotations deal with the following equation-constrained optimization problem:
\begin{align*}
  \underset{T \in \mathbb{R}^{k \times k}}{\mbox{argmax}} \; Q_{\omega}(AT) \mbox{ subject to } T^{\top} T = I_k,
\end{align*}
where $I_k$ is the $k$-th identity matrix, and $A \in \mathbb{R}^{p \times k}$ is an initial solution estimated via methods such as maximum likelihood or least squares. The matrix $\Lambda$ is defined as $\Lambda = (\lambda_{ij})^{1\leq i \leq p}_{1 \leq j \leq k} = AT$, and the function $Q_{\omega}(AT) = Q_{\omega}(\Lambda)$ is defined as 
\begin{align*}
  Q_{\omega}(\Lambda)
  =
  \sum_{i=1}^p \sum_{j=1}^k \lambda_{ij}^4
  -
  \frac{\omega}{p}
  \sum_{j=1}\left(\sum_{i=1}^p \lambda_{ij}^2 \right)^2,
\end{align*}
for a given hyperparameter $0 \leq \omega \leq p$. Here, $p$ denotes the number of observed variables, and $k$ is the number of common factors. For example, $Q_{0}$ is the quartimax criterion, $Q_{1}$ varimax criterion, $Q_{k/2}$ equamax criterion, and $Q_{(p(k-1))/(p+k-2)}$ parsimax criterion, as summarized in \citep[table 1]{Browne2001}. The main contributions of this study are threefold: one theoretical result and two practical advancements.

The primary contributions of this study are as follows. First, we address theoretical results concerning the existence of orthogonal rotation capable of reconstructing simple structures. Specifically, we establish an equivalent condition for the existence of orthogonal rotation that yield perfect simple structures. Furthermore, an equivalent condition for Thurstone simple structures is analyzed using the theory of polynomial ideals (refer to \ref{AppSec:TSS} for detailed exposition). 

Also, utilizing the equivalent condition for the existence of orthogonal rotation that yield perfect simple structures, we perform Monte Carlo simulations to characterize orthomax criteria systematically. This characterization involves identifying all stationary points of the orthomax criteria, calculating criterion values at these points, and determining the positive or negative definiteness of bordered Hessians. These analyses will allow us to elucidate the functional forms of orthomax criteria, including quartimax, varimax, equamax, and parsimax criteria.

Second, we report an unexpected yet practical finding obtained from our simulation: in orthogonal models, varimax rotation has widespread use due to its accessibility in many statistical software, often as the default or sole option for orthogonal rotation. A factor-loading matrix achieves a perfect simple structure when each row contains at most one nonzero element \citep{BernaardsJennrich2003}. This property is referred to as ``unifactoriality'' by \citep{Kaiser1974} and as ``perfect cluster solution'' by \citep{Browne2001}. Many researchers consider varimax to approximate this property effectively; for instance, \citep[page 113]{Browne2001} states that ``Perfect [simple structure] solutions were handled effectively by varimax.'' However, our Monte Carlo simulations reveal an unanticipated finding, corroborated by algebraic analysis, wherein factor rotations derived using the quartimax criterion exhibit superior interpretability compared to those obtained via the varimax criterion. 

Third, our algebraic approach offers several practical advantages. Unlike numerical approaches, such as the gradient projection method, our algebraic approach enables the computation of all stationary points. This allows not only for the exact solution of the optimization problem for a given criterion and but for the selection of solutions that facilitates interpretability for the analyst. These aspects are elaborated upon in Section 6.

The remainder of this study is structured as follows: Section 2 provides a concise review of orthomax factor rotations. Section 3 establishes an equivalent condition for the existence of orthogonal rotation that yield perfect simple structures. Section 4 extends the discussion to Thurstone simple structures, analyzing the corresponding condition. Section 5 introduces a novel computational algebra algorithm for determining all optimal candidates. Section 6 presents numerical results derived from artificial datasets.

% ------------------------------------------------------------------------------------------------------------------------------- %
 
\section{Orthomax rotations}\label{sec:2}

% ------------------------------------------------------------------------------------------------------------------------------- %

\subsection{Factor rotations}\label{sec:2-1}

Let $\bm{X}=(X_1,\cdots,X_p)^\top$ be a $p$-dimensional random vector. A factor analysis model is expressed as 
\begin{equation*}
\X =\bfmu + \Lambda \F + \eps  , \label{model1}
\end{equation*}
where $\bfmu$ is a mean vector, $\Lambda =(\lambda_{ij})$ is the factor loading matrix, and $\F = (F_1,\cdots,F_k)^\top$ and $\eps  = (\varepsilon_1,\cdots, \varepsilon_p)^\top$ are unobservable random vectors. Here, $A^\top$ denotes the transpose of a matrix $A$. The components of \mbox{\boldmath{$F$}} and \mbox{\boldmath{$\varepsilon$}} are referred to as the {\em common factors} and {\em unique factors}, respectively. It is assumed that the unique factors are mutually uncorrelated, and independent of $\eps$. Additionally, under the orthogonal model assumption, the common factors are uncorrelated.

Let $A \in \mathbb{R}^{p \times k}$ denote the estimated factor loading matrix under the orthogonal model, referred to as the {\em initial solution}. A factor loading matrix is known to exhibit rotational indeterminacy, as both $A$ and $AT$ yield the same covariance matrix, $\Sigma = \mathrm{Var}(\bm{X})$, for any regular matrix $T$.

Hence, factor analysis includes factor rotations, where the aim is to find a regular matrix $T$ that transforms $A$ into a simplified factor loading matrix $\Lambda = AT$. Notable simplicity criteria include the {\em perfect simple structure} and {\em Thurstone simple structure}.

A factor-loading matrix $\Lambda$ has a {\em perfect simple structure} if each row of $\Lambda$ contains at most one nonzero element, as defined by \citep{BernaardsJennrich2003}. This structure is alternatively referred to as {\em unifactoriality} \citep{Kaiser1974} or {\em a perfect cluster solution} \citep{Browne2001}.

In addition, {\em Thurstone's rules} \citep{thurstone1947} for the simple structure of a factor matrix $\Lambda$ with $k$ columns, as described by \citep[section ``Simplicity of a Factor Pattern'']{Browne2001} are as follows:
\begin{enumerate}
\item\label{enu:Thurstone1}
  each row of $\Lambda$ contains at least one zero element, 
\item\label{enu:Thurstone2}
  each column of $\Lambda$ contains at least $k$ zero elements,
\item\label{enu:Thurstone3}
  each pair of columns of $\Lambda$ has several rows with a zero element in one column and a nonzero in the other, 
\item\label{enu:Thurstone4}
  for $k \geq 4$, every pair of columns of $\Lambda$ has several rows with zero elements in both columns, and 
\item\label{enu:Thurstone5}
  each pair of columns of $\Lambda$ has few rows with nonzero elements in both columns.
\end{enumerate}
A factor matrix $\Lambda$ is said to have a {\em Thurstone simple structure} if thse conditions are satisfied.

To derive $\Lambda = AT$ with the specified simple structures for a given initial solution $A$, numerous rotation methods have been proposed. If $T$ is orthogonal, the process is referred to as an {\em orthogonal rotation}. This study emphasizes orthogonal rotations and their corresponding criteria.

% ------------------------------------------------------------------------------------------------------------------------------- %

\subsection{Orthogonal rotations}\label{sec:2-2}

Let $Q$ represent an orthogonal criterion. To derive $\Lambda = AT$ with simple structures for a given initial solution $A$, the following equation-constrained optimization problem is solved:
\begin{align}
  \underset{T \in \mathbb{R}^{k \times k}}{\mbox{argmax}} \; Q(AT) \mbox{ subject to } T^{\top} T = I_k,
  \label{opt:ortho}
\end{align}
where $Q(AT)$ denotes the orthogonal rotation criterion evaluated at $AT$. We note that $T$ is orthogonal if and only if $T^{\top} T = I_k$. Some standard rotations are formulated as minimization problems. When necessary we will re-formulate them as equivalent maximization problems.

To maximize $Q$ over orthogonal matrices $T \in \mathbb{R}^{k \times k}$, the existing literature widely adopts the gradient projection algorithm described in \citep{Jennrich2001}. The stopping condition for this algorithm derives directly from \citep[equation (7)]{Jennrich2001} and will be utilized in Section 5. We conclude this section with the stopping condition.
\begin{pro}\label{thm:eq7}
  The orthognal rotation criterion $Q$ has a stationary point at $T \in \mathbb{R}^{k \times k}$ restricted to $\{ T \in \mathbb{R}^{k \times k} : T^{\top} T = I_k\}$ if and only if $T^{\top} \frac{\partial Q(AT)}{\partial T}$ is a symmetric matrix, where
  \begin{align}\label{eq:algeq-pro}
    \frac{\partial Q(AT)}{\partial T} = \left(\frac{\partial Q(AT)}{\partial t_{jl}}\right)_{\substack{1 \leq j \leq k \\ 1 \leq l \leq k}}
  \end{align}
  denotes the gradient of $Q$ at $T = (t_{jl})^{1 \leq j \leq k}_{1 \leq l \leq k}$;
\end{pro}

% ------------------------------------------------------------------------------------------------------------------------------- %

\subsection{Orthomax criteria}\label{sec:2-3}

Various orthogonal rotation criteria can be expressed by the {\em orthomax criterion} $Q_{\omega}(\Lambda)$, which is defined as
\begin{align}\label{eq:orthomax}
  Q_{\omega}(\Lambda)
  =
  \sum_{i=1}^p \sum_{j=1}^k \lambda_{ij}^4
  -
  \frac{\omega}{p}
  \sum_{j=1}^k\left(\sum_{i=1}^p \lambda_{ij}^2 \right)^2,
\end{align}
where $\Lambda = (\lambda_{ij})^{1 \leq i \leq p}_{1 \leq j \leq k}= AT$. For instance, $Q_{0}$ is the quartimax criterion, $Q_{1}$ is the varimax criterion, $Q_{k/2}$ is the equamax criterion, and $Q_{(p(k-1))/(p+k-2)}$ is the parsimax criterion as shown in \cite[table 1]{Browne2001}. Note that $\kappa$ in \cite[table 1]{Browne2001} corresponds to $\omega/p$ in \eqref{eq:orthomax}. Given $A$ is a matrix with real coefficients, the orthomax criterion $Q_{\omega}(\Lambda) = Q_{\omega}(AT)$ can be expressed as a polynomial in the indeterminates $t_{jl}$, with real coefficients for $1 \leq j, l \leq k$, where $T = (t_{jl})^{1 \leq j \leq k}_{1 \leq l \leq k}$. Hence,  
\begin{align*}
  Q_{\omega}(\Lambda) \in \mathbb{R}[t_{jl} : 1 \leq j, l \leq k].
\end{align*}
Here, $\mathbb{R}[t_{jl} : 1 \leq j, l \leq k]$ denotes the polynomial rings of indeterminates $t_{jl}$ with real coefficients (see \ref{sec:a1} for details). Consequently, Proposition \ref{thm:eq7} provides a system of algebraic equations applicable to any orthomax criterion $Q_{\omega}(\Lambda)$, which is an equivalent condition for stationary points. In other words, the proposition provides a system that allows algebraic operations for any orthomax criterion to find all stationary points. Thus, this proposition facilitates computation of all solutions to the system, that is stationary points, for any orthomax criterion $Q_{\omega}(\Lambda)$ using the algebraic operations within $\mathbb{R}[t_{jl} : 1 \leq j, l \leq k]$. Hence, all candidates for the optimization problem \eqref{opt:ortho} can be derived. In particular, since the space $\{ T \in \mathbb{R}^{k \times k} : T^{\top} T = I_k\}$ is compact, we can seek a candidate to maximize $Q_{\omega}(\Lambda)$.

We conclude this section by presenting the property concerning perfect simple structures and orthomax criteria, as established in \cite[theorem 1]{BernaardsJennrich2003}.
\begin{pro}\label{thm:perfect-simple}
  Consider any orthomax criterion $Q_{\omega}$ with $0 \leq \omega \leq p$. If $T$ is an orthogonal matrix and $\Lambda = AT$ has a perfect simple structure, then $\Lambda$ maximizes the criterion $Q_{\omega}$ over all the orthogonal matrices. Furthermore, if $A$ has full column rank, any rotation of $A$ that maximizes the criterion differs from $\Lambda$ by at most one column permutation, and the column sign changes. 
\end{pro}

% ------------------------------------------------------------------------------------------------------------------------------- %
% ------------------------------------------------------------------------------------------------------------------------------- %
 
\section{Perfect simple structure}\label{sec:3}

Proposition \ref{thm:perfect-simple} describes the special properties of perfect simple structures. Therefore, in this section, we discuss the properties of perfect simple structures. Monte Carlo simulations are conducted in Section \ref{sec:6} based on the properties presented in this section.

We present the following theorem. Proposition \ref{thm:perfect-simple} assumes that there exists an orthogonal matrix $T$ such that $\Lambda = AT$ has a perfect simple structure. The following theorem provides the equivalent conditions. 
\begin{thm}\label{thm:perfect-simple-iff}
  The following equations are equivalent:
  \begin{enumerate}
  \item\label{thm:perfect-simple-iff-1}
    there exists an orthogonal matrix $T$ such that $\Lambda = AT$ has a perfect simple structure,
  \item\label{thm:perfect-simple-iff-2}
    there exist $k$ or fewer clusters consisting of rows of $A$ such that 
    \begin{enumerate}
    \item\label{thm:perfect-simple-iff-2-1}
      rows in the same cluster are parallel to each other, and 
    \item\label{thm:perfect-simple-iff-2-2}
      rows in different clusters are orthogonal to each other.
    \end{enumerate}
  \end{enumerate}
  Here, an element of a partition of rows is referred to as a cluster.
  \begin{proof}
    First, suppose that there exists an orthogonal matrix $T$ such that $\Lambda = AT$ has a perfect simple structure. Then, each row of $\Lambda$ has at most one non-zero element. Thus, there exist $k$ or fewer clusters consisting of rows of $\Lambda$ that satisfy Conditions 2(a) and 2(b). Because $T$ is orthogonal, $T^{-1} = T^{\top}$ is also orthogonal. Orthogonal matrices do not affect inner products. Therefore, there exist $k$ or fewer clusters consisting of rows of $A = \Lambda T^{\top}$ that satisfy conditions 2(a) and 2(b).

    Second, suppose that there exist $k$ or fewer clusters $c_1, \ldots, c_m$ of $A$ that satisfy conditions 2(a) and 2(b), where $m \leq k$ is the number of clusters. Moreover there exist row vectors $c_{m+1}, \ldots, c_k \in \mathbb{R}^k$ such that $c_j$ and $c_l$ are orthogonal to each other for different $1 \leq j,l \leq k$, according to the basis extension theorem and the Gram-Schmidt process. Let $s_j = c_j/|c_j|$. We construct an orthogonal matrix as like
    \begin{align*}
      S
      =
      \begin{pmatrix}
        s_1
        \\
        \vdots
        \\
        s_k
      \end{pmatrix}.
    \end{align*}
    Let $a_i$ be the $i$-th row of $A$ for each $1 \leq i \leq p$. Suppose that $a_i$ is parallel to the $j$-th cluster $c_j$. The angle $\theta$ between $a_i$ and $c_j$ is either $0$ or $\pi$. As $a_i$ and $s_l$ are orthogonal for each $l \neq j$, the inner products $a_i s_l^{\top} = 0$. Therefore
    \begin{align}
      a_i S^{\top}
      &=
      a_i
      \begin{pmatrix}
        s_1^{\top} & \cdots & s_k^{\top}
      \end{pmatrix}
      \notag
      \\
      &=
      \begin{pmatrix}
        a_i s_1^{\top} & \cdots & a_i s_k^{\top}
      \end{pmatrix}
      \notag
      \\
      &=
      \begin{pmatrix}
        0 & \cdots & 0 & \overbrace{a_i s_j^{\top}}^{j\mbox{-th}} & 0 & \cdots & 0
      \end{pmatrix}
      \label{eq:rep-sim}
      \\
      &=
      \begin{cases}
        \begin{pmatrix}
          0 & \cdots & 0 & \overbrace{|a_i|}^{j\mbox{-th}} & 0 & \cdots & 0
        \end{pmatrix}
        &
        (\theta = 0)
        \\
        \begin{pmatrix}
          0 & \cdots & 0 & \overbrace{-|a_i|}^{j\mbox{-th}} & 0 & \cdots & 0
        \end{pmatrix}
        &
        (\theta = \pi)
      \end{cases}
      \notag
    \end{align}
    holds. 
    Thus, we obtain the orthogonal matrix $T = S^{\top}$ such that $\Lambda = AT$ has a perfect simple structure.
  \end{proof}
\end{thm}

%Next, we consider the case in which the global optimum is not uniquely determined for orthogonal rotations, based on Theorem \ref{thm:perfect-simple-iff}.

%\begin{cor}\label{cor:perfect-simple-iff}
%    If the number $m$ of clusters consisting of rows of $A$ satisfying the conditions 2(a) and 2(b) is less than the number of common factors $k$, there exist an infinite number of orthogonal matrices $T$ such that $\Lambda = AT$ has a perfect simple structure.
%    \begin{proof}
%      Assume $m < k$. Since there exists an orthogonal matrix such that the column vectors of $AT$ from the first to the $m$th are nonzero vectors, and the column vectors of $AT$ from the $(m+1)$th to the $k$th are zero vectors, as in Eq. \eqref{eq:rep-sim} of the proof of Theorem \ref{thm:perfect-simple-iff}, let $T$ be such an orthogonal matrix. For $m+1 \leq j \leq k$, we focus on the $j$-th column of $AT$. Then, we have a system of linear equations:
%      \begin{align*}
%        A
%        \begin{pmatrix}
%          t_{1j}
%          \\
%          \vdots
%          \\
%          t_{pj}
%        \end{pmatrix}
%        =
%        \begin{pmatrix}
%          0
%          \\
%          \vdots
%          \\
%          0
%        \end{pmatrix}.
%      \end{align*}
%      Since we have $\mathrm{rank}(A) = m < k$ by the assumption, the solutions of each column vector $(t_{1j},\ldots,t_{pj})^{\top}$ to the above system of linear equations can be represented $k-m$ mediator variables. Thus, an infinite number of orthogonal matrices exist, such that $\Lambda = AT$ has a perfect simple structure.
%  \end{proof}
%\end{cor}

We conclude this section with the following remark.

\begin{rem}
  As the second claim of Proposition \ref{thm:perfect-simple} indicates, if $k=m$ is satisfied in Theorem \ref{thm:perfect-simple-iff}, orthogonal matrices such that $\Lambda = AT$ has a perfect simple structure that differs by at most column permutations and column sign changes. In fact, the orthogonal matrix $S$ constructed in the proof of Theorem \ref{thm:perfect-simple-iff} is uniquely determined, except for column permutations and column sign changes when $k=m$.
\end{rem}

This section reveals the form of the initial solutions that make perfect simple structures the global optima. Since we can determine the optima except for column permutations and column sign changes in the case $k=m$, we provide an initial solution that satisfies conditions 2(a), 2(b) and $k = m$ for our Monte Carlo simulations in Section \ref{sec:6}.

% ------------------------------------------------------------------------------------------------------------------------------- %
% ------------------------------------------------------------------------------------------------------------------------------- %
  
\section{Thurstone simple structure}\label{sec:4}

In this section, we characterize Thurstone simple structures to assess their optimality concerning the orthomax criteria. Prior to this characterization, we clarify ambiguities in Thurstone's rules \ref{enu:Thurstone3} and \ref{enu:Thurstone4} by introducing lower limits, denoted as $\gamma$ and $\delta$:
\begin{enumerate}
\setcounter{enumi}{2}
\item\label{enu:Thurstone3-a}
  For every pair of columns of $\Lambda$, at least $\gamma$ rows must contain a zero element in one column but not in the other, 
\item\label{enu:Thurstone4-a}
  for $m \geq 4$, every pair of columns of $\Lambda$ must have at least $\delta$ rows with zero elements in both columns.
\end{enumerate}
Since the satisfaction of rules \ref{enu:Thurstone3} and \ref{enu:Thurstone4} implies compliance with rule \ref{enu:Thurstone5}, we do not address ambiguities in the rule \ref{enu:Thurstone5}. For this study, we state that $\Lambda$ has a {\em Thurstone simple structure of class $(\gamma,\delta)$} if $\Lambda$ satisfies rules \ref{enu:Thurstone1}, \ref{enu:Thurstone2}, \ref{enu:Thurstone3-a} and \ref{enu:Thurstone4-a}. Note that $\Lambda$ has a perfect simple structure in the case such that $\Lambda$ has a Thurstone simple structure of class $(\gamma,\delta)$ such that $\gamma+\delta=p$.

The following theorem, stated in \ref{AppSec:TSS}, establishes that Thurstone simple structures lack the special properties inherent to perfect simple structures.

\begin{thm}\label{thm:thurstone-simple}
  Let $Q_{\omega}$ be any orthomax criterion with $0 \leq \omega \leq p$, and let $1 \leq \gamma \leq p$ and $1 \leq \delta \leq p$. If $T$ is an orthogonal matrix, and even if $\Lambda = AT$ has Thurstone simple structure of class $(\gamma,\delta)$, then $\Lambda$ does not maximize the criterion $Q_{\omega}$ over all orthogonal matrices unless additional conditions are imposed when $\gamma + \delta \neq p$.
\end{thm}

If $\gamma + \delta = p$ holds, as outlined in the subsequent remark, the Thurstone simple structure of class $(\gamma,\delta)$ is attained as the global optimum.

\begin{rem}
  If $\gamma + \delta = p$ holds in Theorem \ref{thm:thurstone-simple}, then $\Lambda$ maximizes the criterion $Q_{\omega}$ among all orthogonal matrices. Notably, $\Lambda$ has a perfect simple structure when it satisfies the Thurstone simple structure of class $(\gamma,\delta)$ with $\gamma + \delta = p$, as established in Proposition \ref{thm:perfect-simple}.
\end{rem}

Section \ref{sec:6} presents Monte Carlo simulations to analyze the behavior of global optima and stationary points under incremental deviations from a perfect simple structure. Therefore, in the next section, we design an algorithm that is capable of computing not only global optima but all stationary points for the orthomax criterion.

% ------------------------------------------------------------------------------------------------------------------------------- %
% ------------------------------------------------------------------------------------------------------------------------------- %
 
\section{Algebraic approach for Orthomax rotations}\label{sec:5}

To characterize the orthomax criteria $Q_{\omega}(\Lambda)$ based on global optima and all stationary points for the optimization problem \eqref{opt:ortho}, which maximizes $Q_{\omega}(\Lambda)$, we develop two algorithms in this section. Numerical approaches may yield local optima that are not global optima. Therefore, an algebraic approach is adopted. Prior to presenting the algorithm, we establish the following corollary derived directly from Proposition \ref{thm:eq7}:

\begin{cor}\label{cor:eq7}
  If an orthogonal matrix $T$ is an optimal solution to \eqref{opt:ortho} for an orthogonal rotation criterion $Q$, then $T^{\top} \frac{\partial Q(AT)}{\partial T}$ is symmetric. In other words, if $T$ is the optimum of \eqref{opt:ortho}, $T$ satisfies 
  \begin{align}\label{eq:algeq}
    \begin{cases}
      \left(T^{\top} \frac{\partial Q(AT)}{\partial T}\right)_{jl}
      =
      \left(T^{\top} \frac{\partial Q(AT)}{\partial T}\right)_{lj}
      &
      (1 \leq j < l \leq k)
      \\
      T^{\top} T = I_k
    \end{cases}.
  \end{align}
\end{cor}

Since $\{ T \in \mathbb{R}^{k \times k} : T^{\top} T = I_k\}$ is a compact space, \eqref{opt:ortho} admits optimal solution. In general, local optima may not coincide with global optima. By computing all algebraic solutions to Eq. \eqref{eq:algeq} for $Q = Q_{\omega}$ the global optimum of \eqref{opt:ortho} for $Q = Q_{\omega}$ can be identified among them. Now, we propose Algorithm \ref{alg:glo} to consistently locate global optima.

\begin{algorithm}[h]
\caption{An algorithm to find a global optimum}
    \label{alg:glo}
    \begin{algorithmic}[1]    
    \REQUIRE an initial solution $A$, a hyper parameter $\omega$
    \ENSURE
      global optimum in \eqref{opt:ortho} for $Q = Q_{\omega}$
      and
      stationary points restricted to $\{ T \in \mathbb{R}^{k \times k} : T^{\top} T = I_k\}$ for $Q = Q_{\omega}$
    \STATE\label{alg:glo-1}
      $\mathcal{T} = \mathfrak{T}_1, \ldots, \mathfrak{T}_L$ $\leftarrow$ connected components of all solutions to \eqref{eq:algeq}
%      \begin{align*}\label{eq:algeq-alg}
%        \begin{cases}
%          \left(T^{\top} \frac{\partial Q_{\omega}(AT)}{\partial T}\right)_{jl}
%          =
%          \left(T^{\top} \frac{\partial Q_{\omega}(AT)}{\partial T}\right)_{lj}
%          &
%          (1 \leq j < l \leq k)
%          \\
%          T^{\top} T = I_k
%        \end{cases}
%      \end{align*}
    \FOR{$a = 1, \ldots, L$}\label{alg:glo-2}
      \IF{$\mathfrak{T}_a \subset \mathbb{R}^{k \times k}$ is a component consisting of only one point}\label{alg:glo-3}
        \STATE\label{alg:glo-4}
        $\mathfrak{t}_a$ $\leftarrow$ the element of the connected space $\mathfrak{T}_a \subset \mathbb{R}^{k \times k}$ 
        \STATE\label{alg:glo-5}
        $q_a$ $\leftarrow$ the value $Q_{\omega}(A\mathfrak{t}_a)$
      \ELSE
        \STATE\label{alg:glo-6}
        $\mathfrak{t}_a$ $\leftarrow$ a sample point of the connected space $\mathfrak{T}_a \subset \mathbb{R}^{k \times k}$
        \STATE\label{alg:glo-7} $q_a$ $\leftarrow$ the value $Q_{\omega}(A\mathfrak{t}_a)$
      \ENDIF
    \ENDFOR
    \STATE\label{alg:glo-8} $\ell$ $\leftarrow$ the index that takes the maximum value among $q_1, \ldots, q_L$
    \RETURN\label{alg:glo-9} $(\mathfrak{t}_\ell, \{ \mathfrak{t}_1, \ldots, \mathfrak{t}_L\})$
    \end{algorithmic}
\end{algorithm}

Step \ref{alg:glo-1} involves computing the connected components of all algebraic solutions to Eq. \eqref{eq:algeq} using an algebraic approach. If $\mathfrak{T}_a \subset \mathbb{R}^{k \times k}$ is a singleton component, the orthomax criterion value $q_a$ is evaluated at the element of the connected space in $\mathfrak{T}_a \subset \mathbb{R}^{k \times k}$ as described in Steps \ref{alg:glo-4} and \ref{alg:glo-5}. Otherwise, a sample point $\mathfrak{t}_a$ is selected from the connected space $\mathfrak{T}_a \subset \mathbb{R}^{k \times k}$ and $q_a$ is evaluated at $\mathfrak{t}_a$ as in Steps \ref{alg:glo-6} and \ref{alg:glo-7}. Importantly, $Q_{\omega}(A\mathfrak{t})$ remains invariant for any $\mathfrak{t} \in \mathfrak{T}_a$, since every $\mathfrak{t} \in \mathfrak{T}_a$ is a stationary point restricted to $\{ T \in \mathbb{R}^{k \times k} : T^{\top} T = I_k\}$ of $Q = Q_{\omega}$. At Step \ref{alg:glo-8}, the index $\ell$ corresponding to the maximum $q_1, \ldots, q_L$ is selected. Given that $\{ T \in \mathbb{R}^{k \times k} : T^{\top} T = I_k\}$ is a compact space, the global optimum $\mathfrak{t}_{\ell}$ is obtained as detailed in Step \ref{alg:glo-9}.

Notably, global optima and stationary points of orthomax criteria may not exist as discrete points; they might manifest as curves, surfaces, etc. Algorithm \ref{alg:glo} accommodates this by considering connected components. In the next section, we apply Algorithm \ref{alg:glo} to our Monte Carlo simulations, which confirm that all global optima and stationary points are indeed discrete points.

It is important to note that Eq. \eqref{eq:algeq} establishes an equivalence not for local maxima but for stationary points. Accordingly, we present an algorithm to classify the stationary points \(\mathfrak{t}_1, \ldots, \mathfrak{t}_L\) using the second-order sufficient optimality conditions derived from the bordered Hessian. Specifically, a point satisfying the bordered Hessian criteria for a local minimum is termed a \textit{second-order sufficient local minimizer}, whereas one fulfilling the criteria for a local maximum is designated a \textit{second-order sufficient local maximizer}. Points where the bordered Hessian fails to provide conclusive evidence — those not meeting second-order sufficient conditions — are classified as \textit{second-order indeterminate points} (which may still correspond to local extrema under higher-order analysis). This classification, based solely on second-order information, facilitates a detailed examination of the shape of the criterion. In constructing this algorithm, we employ the properties of bordered Hessians as outlined in \citep[Sections~3.11 and 7.13]{magnus99} and \citep[Theorems~4 and 5]{GerardDebreu}.

To define the bordered Hessians, consider the Lagrange function for optimization problem \eqref{opt:ortho}:
\begin{align*}
  \Phi(\bm{t}, \bm{\mu}) = Q_{\omega}(AT) + \sum_{j=1}^{\frac{k(k+1)}{2}} \mu_{j} g_{j}.
\end{align*}
Here 
\[
G = \left\{ g_1, \ldots, g_{\frac{k(k+1)}{2}} \right\} = \{ g : g = (T^{\top} T - I_k)_{jl}, 1 \leq j \leq l \leq k\} \subset \mathbb{R}[t_{jl} : 1 \leq j, l \leq k]. 
\]
Additionally, $\mu_{i}$ are the Lagrange multipliers, $\bm{t} = (t_{jl} : 1 \leq j, l \leq k)$, and $\bm{\mu} = \left(\mu_{i} : 1 \leq i \leq k(k+1)/2\right)$.

Let $\bm{\mu}_{a} = \left(\mu_{a,i} : 1 \leq i \leq k(k+1)/2\right)$ be the Lagrange multipliers corresponding to the stationary point $\mathfrak{t}_{a}$ for $1 \leq a \leq L$, satisfying
\[
  \frac{\partial \Phi}{\partial t_{jl}}(\bm{t}_a, \bm{\mu}_a) = \frac{\partial \Phi}{\partial \mu_{i}}(\bm{t}_a, \bm{\mu}_a) = 0, \quad \left(1 \leq a \leq L,\; 1 \leq j,l \leq k,\; 1 \leq i \leq \frac{k(k+1)}{2}\right).
\]
For each $1 \leq a \leq L$, let
\[
  \phi_a(\bm{t}) = \Phi(\bm{t}, \bm{\mu}_a).
\]
We define the Hessian of $\phi_a$ concerning the variables $\{t_{jl}\}$ is given by
\[
  A_a =
  \begin{pmatrix}
    \dfrac{\partial^2 \phi_a}{\partial t_{11} \partial t_{11}} & \cdots & \dfrac{\partial^2 \phi_a}{\partial t_{11} \partial t_{1k}} & \cdots & \dfrac{\partial^2 \phi_a}{\partial t_{11}\partial t_{k1}} & \cdots & \dfrac{\partial^2 \phi_a}{\partial t_{11}\partial t_{kk}} \\[1mm]
    \vdots & \ddots & \vdots &  & \vdots & \ddots & \vdots \\[1mm]
    \dfrac{\partial^2 \phi_a}{\partial t_{kk} \partial t_{11}} & \cdots & \dfrac{\partial^2 \phi_a}{\partial t_{kk} \partial t_{1k}} & \cdots & \dfrac{\partial^2 \phi_a}{\partial t_{kk} \partial t_{k1}} & \cdots & \dfrac{\partial^2 \phi_a}{\partial t_{kk} \partial t_{kk}}
  \end{pmatrix}.
\]
Additionally, the constraint gradient matrix is
\[
  B =
  \begin{pmatrix}
    \dfrac{\partial g_1}{\partial t_{11}} & \cdots & \dfrac{\partial g_1}{\partial t_{1k}} & \cdots & \dfrac{\partial g_1}{\partial t_{k1}} & \cdots & \dfrac{\partial g_1}{\partial t_{kk}} \\[1mm]
    \vdots &  & \vdots &  & \vdots &  & \vdots \\[1mm]
    \dfrac{\partial g_{k(k-1)}}{\partial t_{11}} & \cdots & \dfrac{\partial g_{k(k-1)}}{\partial t_{1k}} & \cdots & \dfrac{\partial g_{k(k-1)}}{\partial t_{k1}} & \cdots & \dfrac{\partial g_{k(k-1)}}{\partial t_{kk}}
  \end{pmatrix},
\]
and let $Z \in \mathbb{R}^{k(k-1) \times k(k-1)}$ denote the zero matrix.

For each integer $b$ satisfying
\[
  \frac{k(k+1)}{2}+1 \leq b \leq k^2,
\]
let $A_a^b$ represent the leading principal $b \times b$ submatrix of $A_a$, and let $B^b$ denote the submatrix of $B$ formed by its first $k(k-1)$ rows and $b$ columns. The $b$-th bordered Hessian of $\phi_a$ is then defined as
\[
  H^b_a =
  \begin{pmatrix}
    A_a^b & (B^b)^{\top} \\[1mm]
    B^b   & Z
  \end{pmatrix}.
\]
To test \emph{definiteness} (i.e., nonsingularity of the Hessian on the tangent space), we require that $\det H^b_a$ is nonzero and satisfies the appropriate sign condition.

Now, we establish a property of bordered Hessians relevant to our next algorithm (refer to \cite[Chapter 7, Section 13]{magnus99} for details). 
\begin{thm}\label{thm:BH}
  For any $\bm{t} \in \mathbb{R}^{k^2}$ such that $\bigwedge_{g \in G} g(\bm{t}) = 0$, we have 
\begin{enumerate}
  \item The point $\bm{t}$ is a \textit{second‐order sufficient local maximizer} if and only if, for every integer 
  \[
    b \in \left\{\frac{k(k+1)}{2}+1, \, \frac{k(k+1)}{2}+2,\, \ldots,\, k^2\right\},
  \]
  the corresponding bordered Hessian satisfies
  \begin{equation}\label{eq:second-order-max}
    (-1)^{b}\, \det H^b_a > 0.
  \end{equation}
  
  \item The point $\bm{t}$ is a \textit{second‐order sufficient local minimizer} if and only if, for every integer 
  \[
    b \in \left\{\frac{k(k+1)}{2}+1,\, \frac{k(k+1)}{2}+2,\, \ldots,\, k^2\right\},
  \]
  the corresponding bordered Hessian satisfies
  \begin{equation}\label{eq:second-order-min}
    (-1)^{\frac{k(k+1)}{2}}\, \det H^b_a > 0.
  \end{equation}
  
  \item The point $\bm{t}$ is a \textit{second‐order indeterminate point} (i.e., it does not satisfy the conditions for a local extremum) if and only if there exists some integer 
  \[
    b \in \left\{\frac{k(k+1)}{2}+1,\, \frac{k(k+1)}{2}+2,\, \ldots,\, k^2\right\}
  \]
  for which neither condition \eqref{eq:second-order-max} nor \eqref{eq:second-order-min} holds.
\end{enumerate}
\end{thm}

Furthermore, (Algorithm \ref{alg:sta}) classifies all stationary points $\mathfrak{t}_1, \ldots, \mathfrak{t}_L$ computed by Algorithm \ref{alg:glo} as follows: second-order sufficient local maxima ``max,''  second-order sufficient local minima ``min,'' and second-order indeterminate point ``indeterminate.'' 

\begin{algorithm}[h]
\caption{An algorithm to classify all stationary points into ``max,'' ``min,'' and ``indeterminate''}
    \label{alg:sta}
    \begin{algorithmic}[1]    
    \REQUIRE all stationary points $\mathfrak{t}_1, \ldots, \mathfrak{t}_L$ computed by Algorithm \ref{alg:glo}
    \ENSURE stationary point patterns $\{ {\tt Pattern} \}$, where ${\tt Pattern}$ is ``max,'' ``min,'' or ``indeterminate''
    \STATE $P$ $\leftarrow$ $\{ \}$
    \FOR{$a = 1, \ldots, L$}
      \IF{any $b=\frac{k(k+1)}{2}+1,\ldots,k^2$ satisfy \eqref{eq:second-order-max}}
        \STATE $P$ $\leftarrow$ $P \cup \{ {\tt max} \}$
      \ELSIF{any $b=\frac{k(k+1)}{2}+1,\ldots,k^2$ satisfy \eqref{eq:second-order-min}}
        \STATE $P$ $\leftarrow$ $P \cup \{ {\tt min} \}$
      \ELSE
          \STATE $P$ $\leftarrow$ $P \cup \{ {\tt indeterminate} \}$
      \ENDIF
    \ENDFOR
    \RETURN $P$
    \end{algorithmic}
\end{algorithm}

Although we focus only on orthomax rotations in this study, the algebraic approach proposed in this section is applicable to arbitrary orthogonal rotations whose criteria are algebraic functions and satisfy the required differentiability assumptions.

% ------------------------------------------------------------------------------------------------------------------------------- %
% ------------------------------------------------------------------------------------------------------------------------------- %
 
\section{Monte Carlo Simulations}\label{sec:6}

Numerical algorithms, including the gradient projection method \citep{Jennrich2001}, depend on initial values and are incapable of identifying all stationary points. In contrast, our algebraic approach, formulated in Algorithm 1, operates independently of initial values and exactly computes all stationary points. Consequently, Algorithm 2, which utilizes Algorithm 1 to classify stationary points, provides a rigorous framework for analyzing the shapes of stationary points in factor rotation criteria. In this section, we present a Monte Carlo simulation employing Algorithms 1 and 2 to examine the characteristics of orthomax criteria. Specifically, we investigate the quartimax criterion $Q_{0}$, varimax criterion $Q_{1}$, equamax criterion $Q_{m/2}$, and parsimax criterion $Q_{(p(m-1))/(p+m-2)}$. Prior to discussing the Monte Carlo simulation results, we detail the initial solutions used in our Monte Carlo simulation. These initial solutions were derived from Theorem \ref{thm:perfect-simple-iff}, as described below.

Firstly, we consider the following initial solution:
\begin{align*}
  A
  =
  \left(
  \begin{array}{rrr}
    0.50 \times 1.0 & 0.40 \times 1.0 & 0.10 \times 1.0
    \\
    0.50 \times 1.1 & 0.40 \times 1.1 & 0.10 \times 1.1
    \\
    0.50 \times 1.2 & 0.40 \times 1.2 & 0.10 \times 1.2
    \\
    0.40 \times 1.0 & - 0.60 \times 1.0 & 0.40 \times 1.0
    \\
    0.40 \times 1.2 & - 0.60 \times 1.2 & 0.40 \times 1.2
    \\
    0.40 \times 0.6 & - 0.60 \times 0.6 & 0.40 \times 0.6
    \\
    0.33 \times 1.0 & - 0.24 \times 1.0 & 0.69 \times 1.0
    \\
    0.33 \times 1.2 & - 0.24 \times 1.2 & 0.69 \times 1.2
    \\
    0.33 \times 1.1 & - 0.24 \times 1.1 & 0.69 \times 1.1
  \end{array}
  \right).
\end{align*}
We systematically vary the entries of the initial solutions to evaluate the properties of global optima, all stationary points, and other related aspects. The initial solution $A$ satisfies condition \ref{thm:perfect-simple-iff-2} of Theorem \ref{thm:perfect-simple-iff}, thereby guaranteeing the existence of an orthogonal matrix $T$ such that $\Lambda = AT$ has a perfect simple structure. Additionally, the number of clusters formed by rows of $A$, which adhere to conditions 2(a) and 2(b), equal the number of common factors $k = 3$. Consequently, the global optimum can be determined up to column permutations and column sign changes.

Next, we consider two methodologies to incrementally disrupt the parallel clusters of the initial solution $A$ in 27 stages. Define
\begin{align*}
  S
  =
  \left(
  \begin{array}{rrr}
1 & 4 & 7 \\
10 & 13 & 16 \\
19 & 22 & 25 \\
2 & 5 & 8 \\
11 & 14 & 17 \\
20 & 23 & 26 \\
3 & 6 & 9 \\
12 & 15 & 18 \\
21 & 24 & 27
  \end{array}
  \right),
  \quad
  W
  =
  \left(
  \begin{array}{rrr}
1 & 2 & 3 \\
4 & 5 & 6 \\
7 & 8 & 9 \\
10 & 11 & 12 \\
13 & 14 & 15 \\
16 & 17 & 18 \\
19 & 20 & 21 \\
22 & 23 & 24 \\
25 & 26 & 27
  \end{array}
  \right).
\end{align*}
Let $U_{ij}$ $(1 \leq i \leq 9, \ 1 \leq j \leq 3)$ be independent and identically distributed (i.i.d) random variables sampled from the uniform distribution $U(-1,1)$. These variables represent perturbation introduced to the initial loading matrix $A$ to induce variability. Then, for a given integer $\ell \ (\ell=1,\dots,27)$, the perturbations $U_{ij}$ are incrementally added to $A_{ij}$ as follows:
$$
(A_\ell^S)_{ij}=A_{ij}+U_{ij}\,\mathbf{1}\{S_{ij}\leq\ell\},\qquad
(A_\ell^W)_{ij}=A_{ij}+U_{ij}\,\mathbf{1}\{W_{ij}\leq\ell\},
$$
where $\mathbf{1}(\cdot)$ denotes the indicator function. Perturbations $U_{ij}$ that do not satisfy 
$$
\sum_{j=1}^3 (A_\ell^S)_{ij}^2 \in [0,1], \quad \sum_{j=1}^3 (A_\ell^W)_{ij}^2 \in [0,1]
$$
are discarded, and new random values are sampled iteratively until these condition are satisfied.

As $\ell$ increases, the entries of $A$ are progressively influenced by noise $U_{ij}$, resulting in a gradual dissolution of the parallel clusters. Specifically, the parallel clusters in $A_{\ell}^{S}$ vanish for $\ell \geq 12$, whereas in $A_{\ell}^{W}$, they vanish for $\ell \geq 22$. Thus, the loss of parallelism in $A_{\ell}^{S}$ occurs more abruptly compared to $A_{\ell}^{W}$. The initial solutions $A_{\ell}^S$ and $A_{\ell}^W$ are referred to as {\tt Type S} and {\tt Type W}, respectively.

We generated 50 sets of $U=(U_{ij})$ to compare results derived using the {\tt GPArotation} package, global optima identified by Algorithm \ref{alg:glo}, and stationary points obtained using Algorithm \ref{alg:glo}. The {\tt GPArotation} package was implemented in {\tt R} and was based on \citep{Jennrich2001}. All parameters, including the threshold for the convergence assessment, were set to their default values. Algorithm \ref{alg:glo} was implemented in a way that allows the computer algebra system {\tt Mathematica} to be executed from the computer algebra system  {\tt SageMath}. In particular, our preliminary investigation revealed that, for initial solutions like those generated above, only a finite number of algebraic solutions existed. Therefore, our implementation does not employ Steps 6, 7, and 8 of Algorithm \ref{alg:glo}. On the other hand, for initial solutions in which the number of clusters satisfying conditions 2(a) and 2(b) of Theorem 1 was less than the number of factors, an infinite number of algebraic solutions were found to exist. Accordingly, when performing such simulations, it is necessary to include Steps 6, 7, and 8 of Algorithm \ref{alg:glo}. The computational time of our algorithm varied from a few minutes to several tens of minutes to obtain results for a single dataset. 

\begin{figure}[!htb]
  \centering
  \includegraphics[scale=0.4]{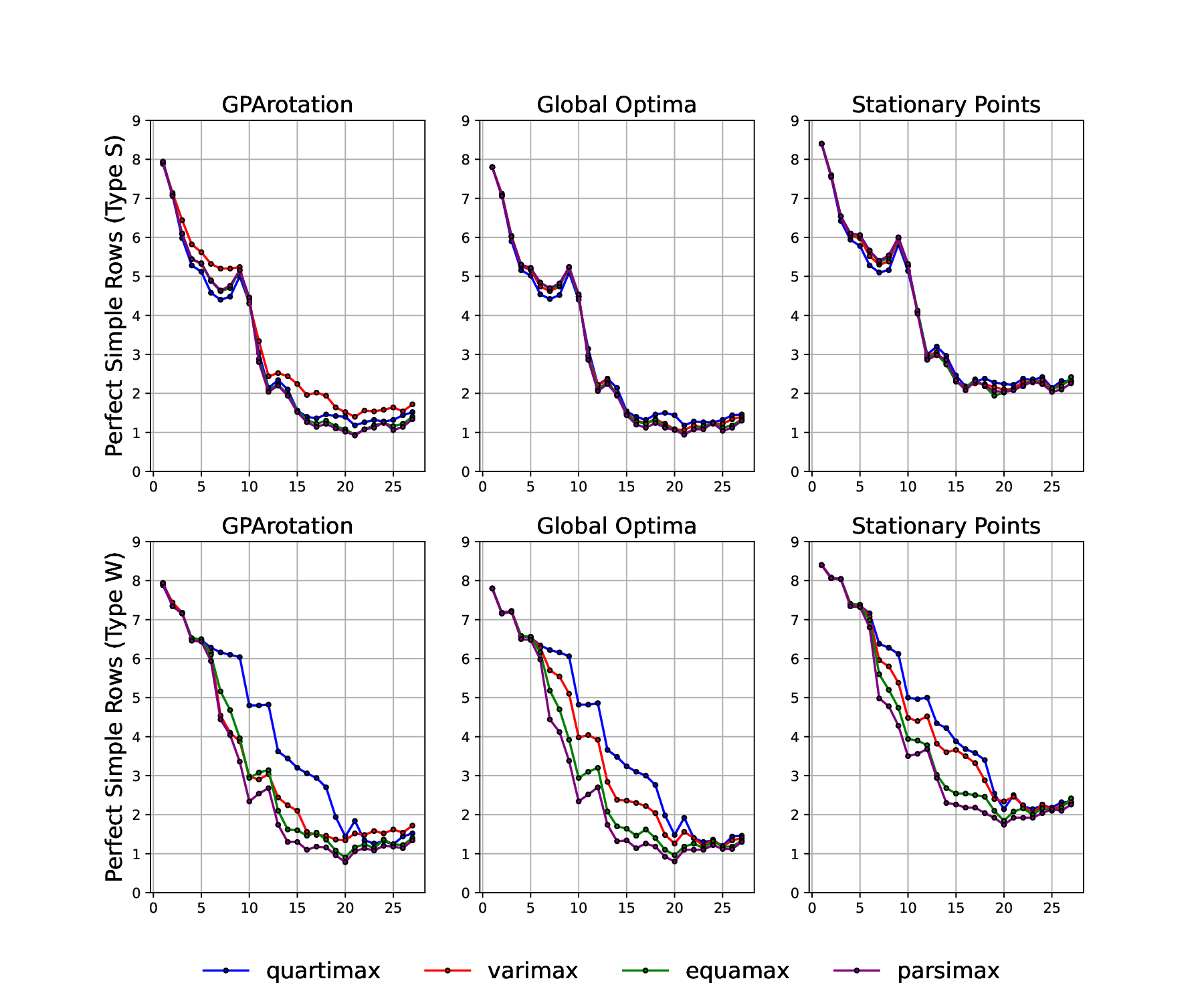}
  \caption{rows such that the absolute values of two or more elements are less than $0.1$}
  \label{fig:perfect-simple-row}
\end{figure}

\begin{figure}[!h]
  \centering
  \includegraphics[scale=0.4]{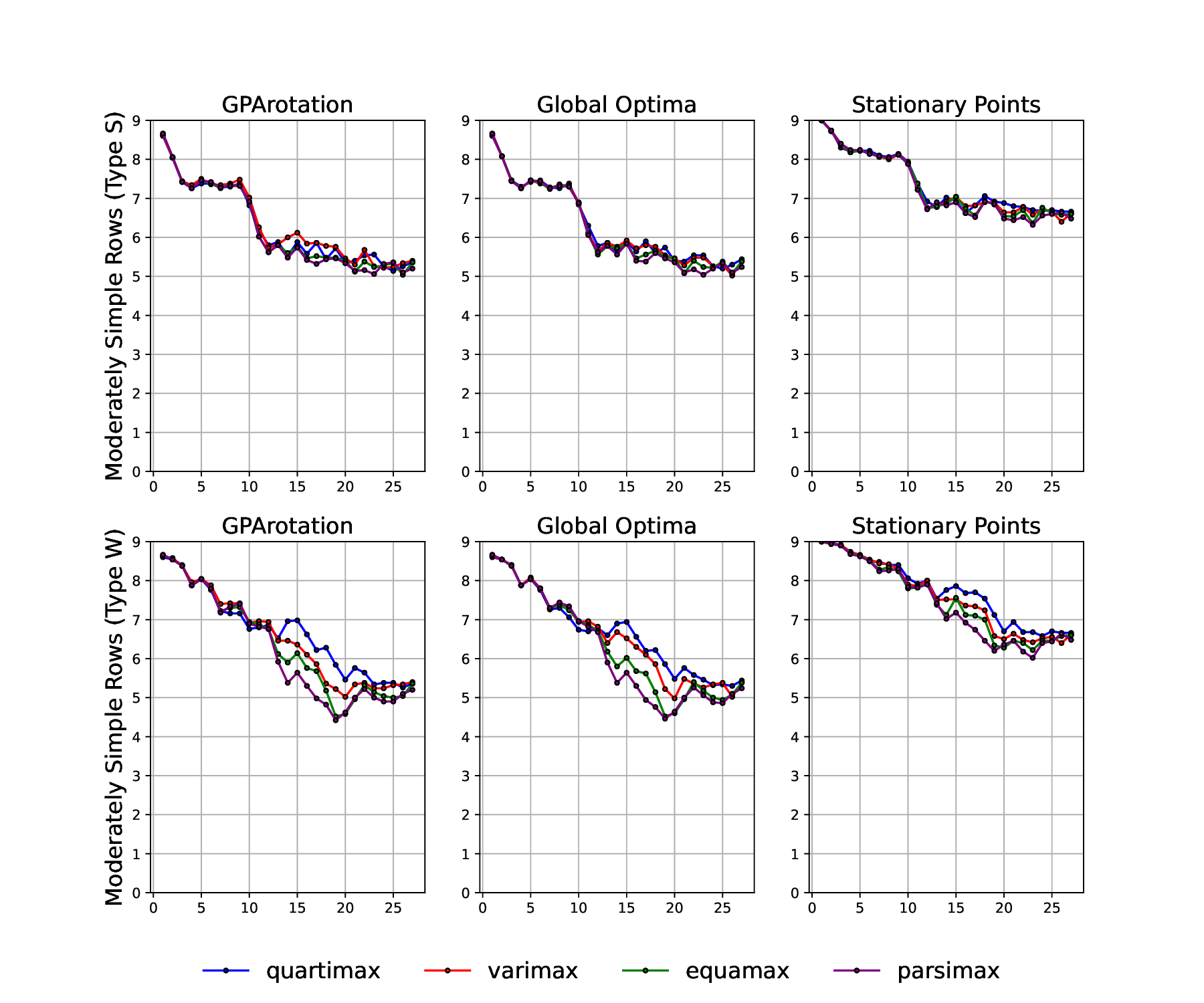}
  \caption{rows such that the absolute values of one or more elements are less than $0.1$.}
  \label{fig:thurstone-simple-row}
\end{figure}

\begin{figure}[!h]
  \centering
  \includegraphics[scale=0.4]{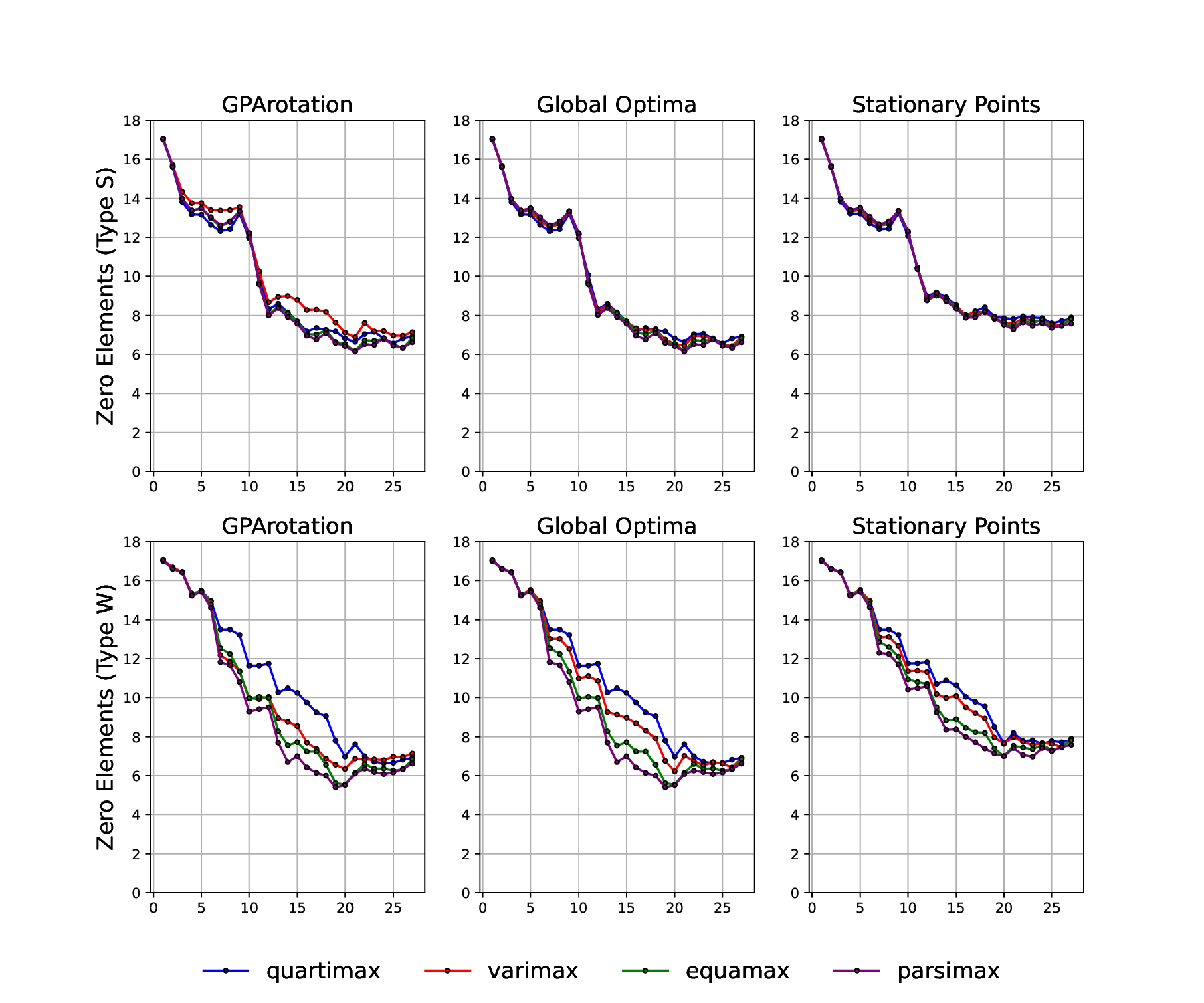}
  \caption{elements whose absolute values are less than $0.1$}
  \label{fig:zero-elements}
\end{figure}

The simulation results are illustrated in Figures \ref{fig:perfect-simple-row} through \ref{fig:zero-elements}.
Figure \ref{fig:perfect-simple-row} depicts the number of rows where the absolute values of at least two elements are less than $0.1$, referred to as ``perfect simple rows.'' Figure \ref{fig:thurstone-simple-row} illustrates the number of rows where the absolute values of at least one element are less than $0.1$, denoted as ``moderately simple rows.'' Figure \ref{fig:zero-elements} displays the number of elements with absolute values less than $0.1$, identified as ``zero elements.'' The horizontal axis in each figure represents the index $\ell \ (\ell = 1, \ldots, 27)$, while the vertical axis indicates the averages of ``perfect simple rows,'' ``moderately simple rows'' and ``zero elements,'' computed across fifty instances of $A_{\ell}^S$ and $A_{\ell}^W$. The figures consists of three panels: the left panels present results from the {\tt GPArotation} package, the middle panels show global optima determined by Algorithm \ref{alg:glo}, and the right panels illustrate results from stationary points maximizing the number of ``perfect simple rows,'' ``moderately simple rows,'' and ``zero elements,'' respectively. The upper panels in each figure correspond to {\tt Type S} results, whereas the lower panels correspond to {\tt Type W} results. Within each panel, the blue lines represent quartimax criterion $Q_{0}$, the red lines varimax criterion $Q_{1}$, the green lines equamax criterion $Q_{m/2}$, and the purple lines parsimax criterion $Q_{(p(m-1))/(p+m-2)}$.

As illustrated in Figures \ref{fig:perfect-simple-row}–\ref{fig:zero-elements}, the behavior of the {\tt GPArotation} output closely resembles that of the global optima computed by Algorithm \ref{alg:glo}. However, {\tt GPArotation} may generally converge to a stationary point that does not correspond to a global optimum. Since Algorithm \ref{alg:glo} is designed to guarantee global optimality, we obtained results that highlight the favorable performance of {\tt GPArotation}.

Moreover, the presence of numerous perfect simple rows, moderately simple rows, and zero elements indicates that the obtained rotation results are highly interpretable. This outcome is unexpected, as the quartimax criterion yields more interpretable rotation results than the varimax criterion, despite the greater popularity of the latter.

Furthermore, a stationary point yields simpler rotation results compared to global optima, particularly in varimax rotation. Notably, these stationary points are not necessarily local maxima. While the orthomax criteria are designed to extract simple structures, they do not fully achieve this objective. Our algebraic framework enables the computation of all stationary points, facilitating the selection of solutions that align with analysts’ interpretability preferences. The following matrices $\Lambda = AT$, present, from left to right: global optima, solutions prioritizing perfect simple rows, solutions focusing on moderately simple rows, and solutions emphasizing zero components. Here, $T$ is selected from the set $t_1, \ldots, t_L$ generated by Algorithm 1. Specifically, the first row corresponds to the an initial {\tt Type W} solution at index $\ell=9$, the second row at $\ell=18$, and the third row at $\ell=27$. Moreover, components with absolute values, truncated to two decimal place, below 0.1 are denoted by \dotuline{~~~~~~~} in the following matrices:

{\tiny\begin{align*}
\begin{array}{l|ccccc}
& \mbox{Global Optima}
& \mbox{Perfect Simple Rows}
& \mbox{Moderately Simple Rows}
& \mbox{Zero Elements}
\\\hline
\ell = 9
&
\begin{pmatrix}
     .12 & -.40 & -.79 \\
    -.12 &  .29 &  .52 \\
     .74 & -.54 & -.79 \\
    -.82 & \dotuline{~~~~~~~} &  \dotuline{~~~~~~~} \\
    -.99 & \dotuline{~~~~~~~} &  \dotuline{~~~~~~~} \\
    -.49 & \dotuline{~~~~~~~} &  \dotuline{~~~~~~~} \\
    \dotuline{~~~~~~~} &  \dotuline{~~~~~~~} & -.79 \\
    \dotuline{~~~~~~~} &  .11 & -.95 \\
    \dotuline{~~~~~~~} &  .10 & -.87 \\
\end{pmatrix}
&
\begin{pmatrix}
     .12 & -.40 & -.79 \\
    -.12 &  .29 &  .52 \\
     .74 & -.54 & -.79 \\
    -.82 & \dotuline{~~~~~~~} &  \dotuline{~~~~~~~} \\
    -.99 & \dotuline{~~~~~~~} &  \dotuline{~~~~~~~} \\
    -.49 & \dotuline{~~~~~~~} &  \dotuline{~~~~~~~} \\
    \dotuline{~~~~~~~} & \dotuline{~~~~~~~} & -.79 \\
    \dotuline{~~~~~~~} &  .11 & -.95 \\
    \dotuline{~~~~~~~} &  .10 & -.87 \\
\end{pmatrix}
&
\begin{pmatrix}
     .12 & -.40 & -.79 \\
    -.12 &  .29 &  .52 \\
     .74 & -.54 & -.79 \\
    -.82 & \dotuline{~~~~~~~} &  \dotuline{~~~~~~~} \\
    -.99 & \dotuline{~~~~~~~} &  \dotuline{~~~~~~~} \\
    -.49 & \dotuline{~~~~~~~} &  \dotuline{~~~~~~~} \\
    \dotuline{~~~~~~~} &  \dotuline{~~~~~~~} & -.79 \\
    \dotuline{~~~~~~~} &  .11 & -.95 \\
    \dotuline{~~~~~~~} &  .10 & -.87 \\
\end{pmatrix}
&
\begin{pmatrix}
     .12 & -.40 & -.79 \\
    -.12 &  .29 &  .52 \\
     .74 & -.54 & -.79 \\
    -.82 & \dotuline{~~~~~~~} &  \dotuline{~~~~~~~} \\
    -.99 & \dotuline{~~~~~~~} &  \dotuline{~~~~~~~} \\
    -.49 & \dotuline{~~~~~~~} &  \dotuline{~~~~~~~} \\
    \dotuline{~~~~~~~} & \dotuline{~~~~~~~} & -.79 \\
    \dotuline{~~~~~~~} &  .11 & -.95 \\
    \dotuline{~~~~~~~} &  .10 & -.87 \\
\end{pmatrix}
\\
\ell = 18
&
\begin{pmatrix}
    -.30 & -.21 & -.82 \\
     .25 &  .14 &  .53 \\
    -.93 & \dotuline{~~~~~~~} & -.82 \\
    -.15 &  .64 &  .11 \\
     .30 & -.35 & -.63 \\
    -.20 & -.40 &  .51 \\
     .14 & \dotuline{~~~~~~~} & -.79 \\
     .17 & \dotuline{~~~~~~~} & -.94 \\
     .15 & \dotuline{~~~~~~~} & -.87 \\
\end{pmatrix}
&
\begin{pmatrix}
    -.79 & -.41 & \dotuline{~~~~~~~} \\
     .52 &  .30 & \dotuline{~~~~~~~} \\
     .33 & -.61 & \dotuline{~~~~~~~} \\
     \dotuline{~~~~~~~} &  .35 &  .56 \\
    -.63 & \dotuline{~~~~~~~} & -.45 \\
     .53 & -.39 & -.15 \\
    -.79 &  .10 & \dotuline{~~~~~~~} \\
    -.95 &  .12 & \dotuline{~~~~~~~} \\
    -.87 &  .11 & \dotuline{~~~~~~~} \\
\end{pmatrix}
&
\begin{pmatrix}
    -.79 & -.41 & \dotuline{~~~~~~~} \\
     .52 &  .30 & \dotuline{~~~~~~~} \\
     .33 & -.61 & \dotuline{~~~~~~~} \\
     \dotuline{~~~~~~~} &  .35 &  .56 \\
    -.63 & \dotuline{~~~~~~~} & -.45 \\
     .53 & -.39 & -.15 \\
    -.79 &  .10 & \dotuline{~~~~~~~} \\
    -.95 &  .12 & \dotuline{~~~~~~~} \\
    -.87 &  .11 & \dotuline{~~~~~~~} \\
\end{pmatrix}
&
\begin{pmatrix}
    -.79 & -.41 & \dotuline{~~~~~~~} \\
     .52 &  .30 & \dotuline{~~~~~~~} \\
     .33 & -.61 & \dotuline{~~~~~~~} \\
     \dotuline{~~~~~~~} &  .35 &  .56 \\
    -.63 & \dotuline{~~~~~~~} & -.45 \\
     .53 & -.39 & -.15 \\
    -.79 &  .10 & \dotuline{~~~~~~~} \\
    -.95 &  .12 & \dotuline{~~~~~~~} \\
    -.87 &  .11 & \dotuline{~~~~~~~} \\
\end{pmatrix}
\\
\ell = 27
&
\begin{pmatrix}
    \dotuline{~~~~~~~} & -.79 &  .42 \\
    \dotuline{~~~~~~~} &  .53 & -.29 \\
     .96 & \dotuline{~~~~~~~} &  .42 \\
     .11 &  .47 &  .46 \\
    -.47 & -.62 &  \dotuline{~~~~~~~} \\
     .41 &  \dotuline{~~~~~~~} & -.53 \\
     .13 & -.15 &  .95 \\
    -.38 &  .44 &  .23 \\
    -.55 &  .56 &  .24 \\
\end{pmatrix}
&
\begin{pmatrix}
     .58 & -.55 & -.40 \\
    -.40 &  .34 &  .30 \\
     .23 &  .67 & -.40 \\
     .35 &  .48 &  .30 \\
     .11 & -.77 & \dotuline{~~~~~~~} \\
    -.50 &  .25 & -.38 \\
     .97 &  \dotuline{~~~~~~~} & \dotuline{~~~~~~~} \\
     \dotuline{~~~~~~~} &  .11 &  .61 \\
     \dotuline{~~~~~~~} &  \dotuline{~~~~~~~} &  .81 \\
\end{pmatrix}
&
\begin{pmatrix}
    -.62 & -.65 & \dotuline{~~~~~~~} \\
     .45 &  .40 & \dotuline{~~~~~~~} \\
    -.71 &  .63 & \dotuline{~~~~~~~} \\
    -.11 &  .39 &  .53 \\
    \dotuline{~~~~~~~} & -.78 & \dotuline{~~~~~~~} \\
    \dotuline{~~~~~~~} &  .36 & -.57 \\
    -.70 & -.12 &  .67 \\
     .34 & \dotuline{~~~~~~~} &  .52 \\
     .50 & \dotuline{~~~~~~~} &  .64 \\
\end{pmatrix}
&
\begin{pmatrix}
    -.62 & -.65 & \dotuline{~~~~~~~} \\
     .45 &  .40 & \dotuline{~~~~~~~} \\
    -.71 &  .63 & \dotuline{~~~~~~~} \\
    -.11 &  .39 &  .53 \\
    \dotuline{~~~~~~~} & -.78 & \dotuline{~~~~~~~} \\
    \dotuline{~~~~~~~} &  .36 & -.57 \\
    -.70 & -.12 &  .67 \\
     .34 & \dotuline{~~~~~~~} &  .52 \\
     .50 & \dotuline{~~~~~~~} &  .64 \\
\end{pmatrix}
\end{array}
\end{align*}}

Thus, utilizing our algebraic approach, various factor rotations can be systematically obtained. This approach enables the computation of all stationary points independently of initial values, presenting a distinct advantage over numerical methods such as the gradient projection method.

\begin{figure}[!ht]
  \centering
  \includegraphics[scale=0.4]{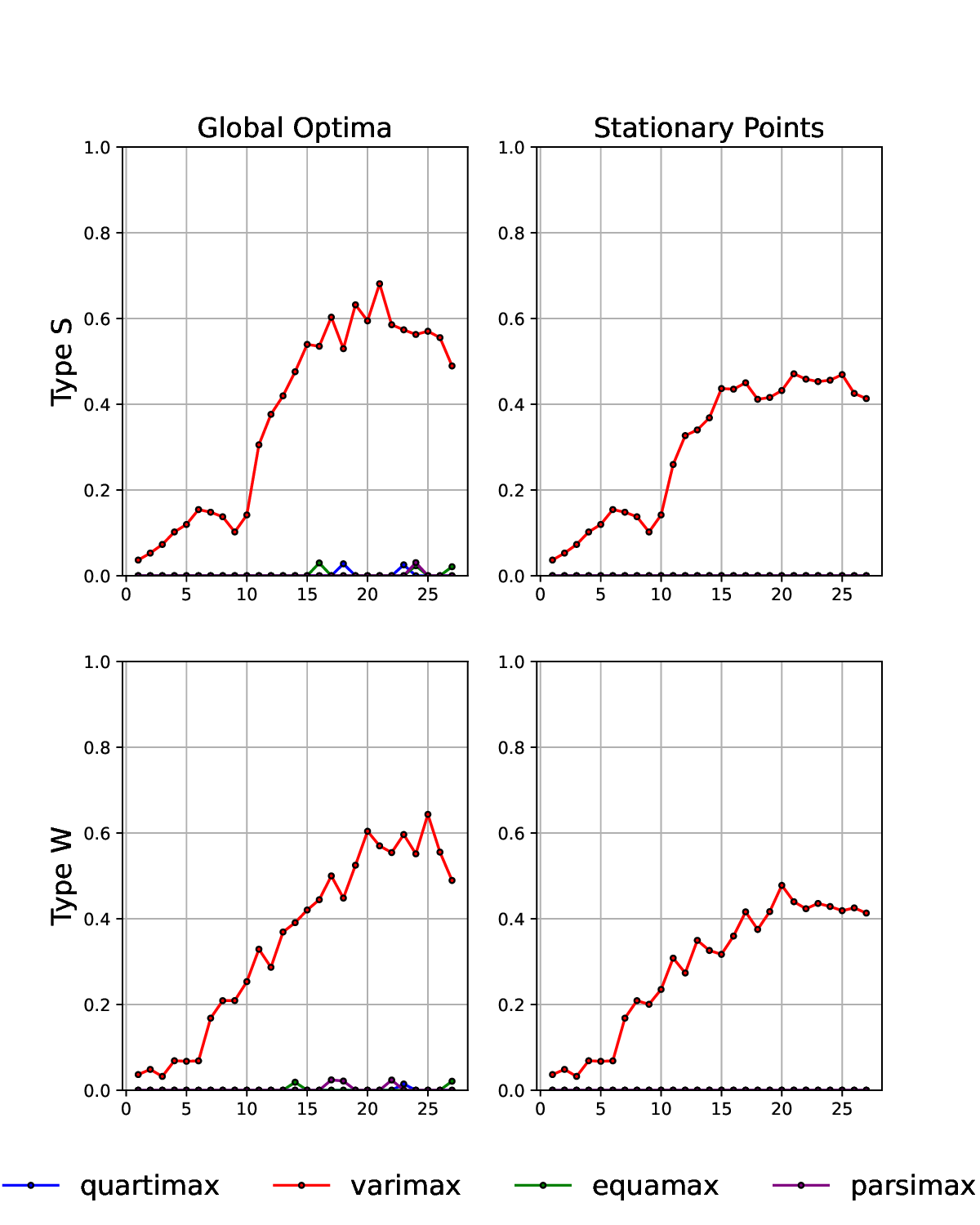}
  \caption{Euclidean distances from the {\tt GPArotation} outputs}
  \label{fig:distance}
\end{figure}

The left panels of Figure \ref{fig:distance} illustrate the mean Euclidean distances between {\tt GPArotation} outputs and global optima, while the right panels display the mean distance between {\tt GPArotation} outputs and the nearest stationary points. Except for the varimax criterion, the Euclidean distance between {\tt GPArotation} outputs and global optima is nearly zero, providing strong evidence that {\tt GPArotation} reliably converges to global optima.

The Euclidean distance corresponding to the varimax criterion may initially appear relatively large in comparison to other criteria. However, the maximum possible Euclidean distance, given by $\sqrt{3 \times 9 \times 2^2} \fallingdotseq 10.392$, indicates that the observed Euclidean distances for the varimax criterion remain relatively small. Moreover, examining the right panel reveals a behavior closely resembling that of the left panel, suggesting that the {\tt GPArotation} has effectively converged to a global optima. Notably, our algebraic approach — specifically, Algorithm \ref{alg:glo} — guarantees the computation of all global optima and stationary points. Consequently, it facilitates the assessment of the performance of existing optimization algorithms, such as the gradient projection method and its implementation in {\tt GPArotation}.

\begin{figure}[!h]
  \centering
  \includegraphics[scale=0.4]{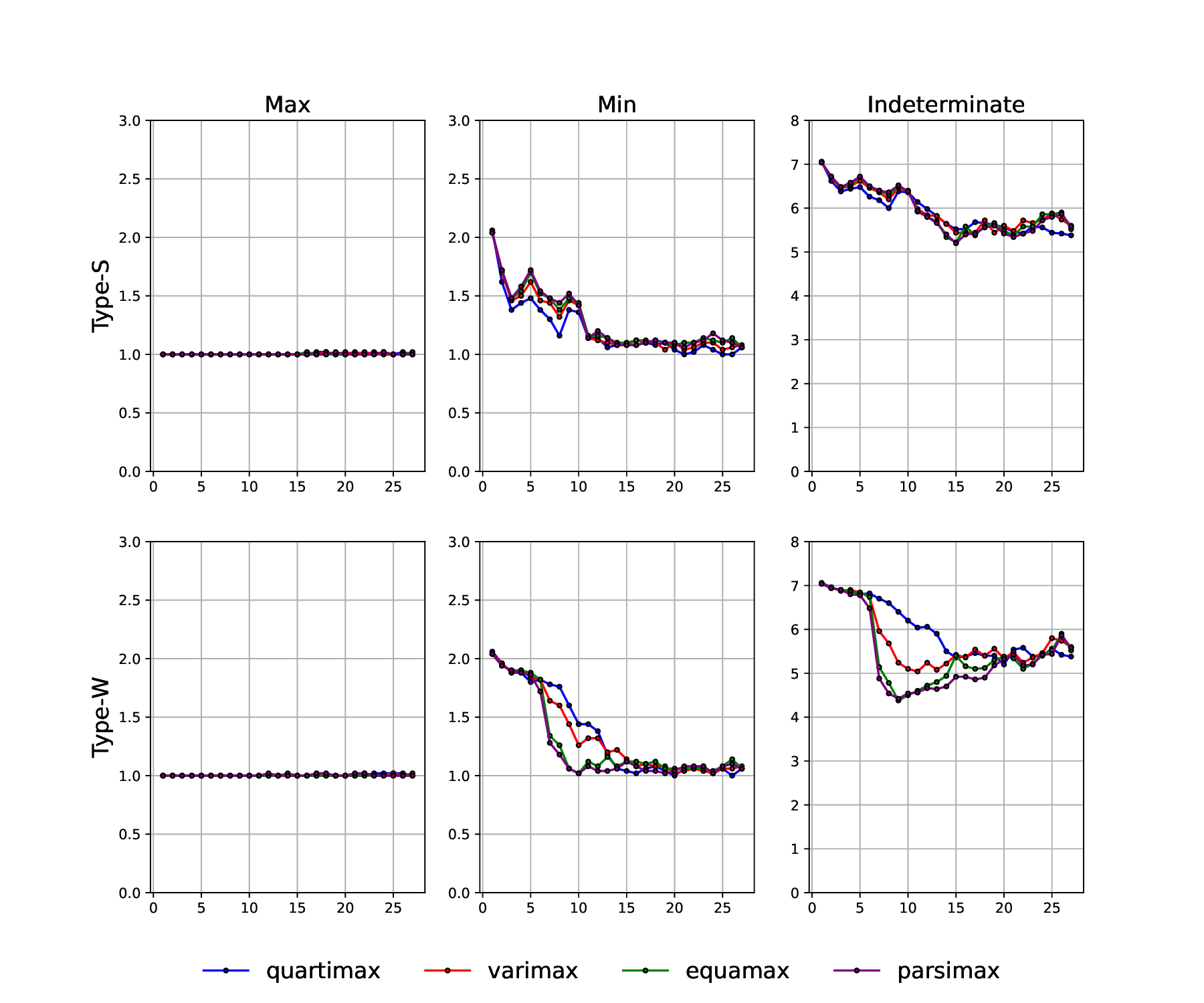}
  \caption{second-order sufficient local maxima ``max,'' second-order sufficient local minima ``min,'' and second-order indeterminate point ``indeterminate''}
  \label{fig:stationary}
\end{figure}

This section concludes with the averages of the numbers of second-order sufficient local maxima ``max,'' second-order sufficient local minima ``min,'' and second-order indeterminate point ``indeterminate'' associated with orthomax criteria. The behavior of orthomax criteria as objective functions is analyzed through these averages, as illustrated in Figure \ref{fig:stationary}.

Orthomax criteria, as objective functions of maximization problems, exhibit a desirable property whereby computational results indicate a unique maximum in most experiments. As demonstrated in Figure \ref{fig:stationary}, most stationary points are classified as ``min'' or ``indeterminate.'' Consequently, even with established algorithms such as the gradient projection method, testing multiple initial values combined with bounded Hessians to ascertain local maximalities, often will guide the algorithm toward global optima. 

% ------------------------------------------------------------------------------------------------------------------------------- %
% ------------------------------------------------------------------------------------------------------------------------------- %
 
\section{Conclusion and future works}\label{sec:7}

In this study, we present the theoretical results concerning perfect simple structures and Thurstone simple structures, formalized in Theorems \ref{thm:perfect-simple-iff} and \ref{thm:thurstone-simple}, respectively. Notably, the Monte Carlo simulation conducted in Section 6 employs initial solutions constructed based on Theorem \ref{thm:perfect-simple-iff}.

Furthermore, we introduce Algorithm 1, which is based on an algebraic approach that is fundamentally distinct from numerical methods, such as the gradient projection technique, which has been extensively utilized. Specifically, numerical approaches exhibit dependency on the initial values, rendering it infeasible to determine all stationary points. In contrast, our algebraic approach operates independently of initial values, thereby ensuring that the output remains unaffected while enabling the computation of all stationary points. Consequently, our algebraic framework provides a systematic method for selecting interpretable solutions from the set of computed stationary points. In particular, Section 6 demonstrates that stationary points having greater simplicity than global optima can be identified, and our algebraic method facilitates factor rotation based on such stationary points.

Moreover, Monte Carlo simulations conducted in this study, reveal that quartimax criterion tends to yield simpler solutions compared to the widely adopted varimax criterion. Additionally, favorable outcomes were observed for the {\tt GPArotation} package, which exhibits a strong tendency to converge to global optima, and for the uniqueness of local maxima across multiple initial solutions in the quartimax, varimax, equamax, and parsimax criteria — an advantageous property in optimization problems involving a maximization objective function.

In this study, we employ an algebraic framework to examine the properties and behavior of orthogonal rotations. Given its broader applicability, future research may extend this approach to analyze oblique rotations, which are more frequently utilized than orthogonal rotations. The insights derived from our investigations of orthogonal rotations are expected to inform such subsequent studies.

Notably, the proposed method demands considerably greater computational resources compared to numerical techniques such as the gradient projection method. Our simulation accounts for three factors; however, as the number of factors increases, computational constraints may arise, rendering solutions infeasible within practical time limits. Therefore, optimizing computational efficiency and memory utilization remains a critical challenge for future research.

\section*{Acknowledgements}
This work was supported by JSPS KAKENHI Grant Number JP23H04474 (K.H.), JP25H01482 (R.F.), JP20K14312 (Y.K.), JP21K18312 (Y.K.), JP23H03352 (M.Y.).

% ------------------------------------------------------------------------------------------------------------------------------- %
% ------------------------------------------------------------------------------------------------------------------------------- %
 
\bibliographystyle{abbrvnat}
\bibliography{main.bib}

% ------------------------------------------------------------------------------------------------------------------------------- %
% ------------------------------------------------------------------------------------------------------------------------------- %
 
\appendix
\def\thesection{Appendix \Alph{section}}

\section{Basic algebraic concepts}\label{sec:a1}

This section presents fundamental concepts in polynomial ring theory, focusing on fields and rings, along with illustrative examples. We begin by formally defining fields.

\begin{dfn}\label{dfn:fields}
A field is a set $F$ equipped with two binary operations ``$+$'' and ``$\cdot$'' defined on $F$ satisfying the following conditions:
\begin{enumerate}
  \item for any $a,b,c \in F$, $(a+b)+c = a+(b+c)$ and $(a \cdot b) \cdot c = a \cdot (b \cdot c)$ (associativity),
  \item for any $a,b,c \in F$, $a \cdot (b + c) = a \cdot b + a \cdot c$ (distributivity),
  \item for any $a, b \in F$, $a + b = b + a$ and $a \cdot b = b \cdot a$ (commutativity),
  \item for any $a \in F$, there exists $0, 1 \in F$ such that $a + 0 = a \cdot 1 = a$ (identities),
  \item for any $a \in F$, there exists $b \in F$ such that $a + b = 0$ (additive inverses),
  \item for any $a \in F$, $a \neq 0$, there exists $b \in F$ such that $a \cdot b = 1$ (multiplicative inverses).
\end{enumerate}
\end{dfn}

For instance, the sets $\mathbb{Q}$, $\mathbb{R}$, and $\mathbb{C}$ are fields, as they satisfy the following conditions with the sum ``$+$'' and product ``$\cdot$''. Conversely, $\mathbb{Z}$ does not form a field as it fails to satisfy the requirement of multiplicative inverses. Indeed, the element $2 \in \mathbb{Z}$ lacks an element $b \in \mathbb{Z}$ such that $2 \cdot b = 1$.

We now proceed to the definition of a commutative ring.

\begin{dfn}
A commutative ring is a set $R = F$ equipped with two binary operations, ``$+$'' and ``$\cdot$'', that satisfy conditions 1-5 outlined in Definition \ref{dfn:fields}.
\end{dfn}

As previously noted, $\mathbb{Z}$ is not a field; however, it is a commutative ring. Furthermore, the set of polynomials is a commutative ring. In this study, we consider algebraic equations of the form \eqref{eq:algeq}, where such equations are defined as $f = 0$, with $f$ being a polynomial whose coefficients belongs to a specified field. So, we provide a formal definition of polynomials.

\begin{dfn}
  A monomial in $\bm{z} = (z_1, \ldots, z_m)$ is an expression of the form $\bm{z}^{\alpha} = z_1^{\alpha_1} \cdots z_m^{\alpha_m}$, where the exponent vector $\alpha = (\alpha_1,\ldots, \alpha_m)$ consists of nonnegative integers, that is $\alpha_i \in \mathbb{Z}_{\geq 0}$. A polynomial $f$ in $\bm{z}$ with coefficients in the real number field $\mathbb{R}$ is a finite linear combination (with coefficients in $\mathbb{R}$) of monomials. Explicitly, we express $f$ as $f = \sum_{\alpha = (\alpha_1,\ldots, \alpha_m) \in \mathbb{Z}_{\geq 0}^m} a_{\alpha} \bm{z}^{\alpha}$, where the summation is taken over a finite set of $m$-tuples $\alpha = (\alpha_1,\ldots, \alpha_m)$. The set of all polynomials in $\bm{z}$ with coefficients in $\mathbb{R}$ is denoted $\mathbb{R}[\bm{z}]$.
\end{dfn}

As noted above, $\mathbb{R}[\bm{z}] = \mathbb{R}[z_1, \ldots, z_m]$ is a commutative ring. Specifically, $\mathbb{R}[\bm{z}]$ is referred to as the polynomial ring. Analogous to the treatment of linear equations via linear subspaces, algebraic equations are systematically addressed through ideals in a polynomial ring. We now proceed to define ideals. In general, ideals are defined for arbitrary ring. For simplicity, we present the definition of ideals for commutative rings.

\begin{dfn}
  Let $R$ be a commutative ring. A subset $\mathcal{I} \subset R$ is an ideal if it satisfies the following conditions:
  \begin{enumerate}
  \item $0 \in \mathcal{I}$,
  \item if $a, b \in \mathcal{I}$, then $a + b \in \mathcal{I}$,
  \item if $a \in \mathcal{I}$ and $b \in R$, then $a \cdot b \in \mathcal{I}$.
  \end{enumerate}
\end{dfn}

The concept of ideals parallels that of linear subspaces, as both structures are closed under addition and multiplication. However, an essential distinction lies in the multiplicative; for a linear subspace, we multiply an element in the field, whereas for ideals, we multiply an element in the ring. 

Just as linear equations are addressed via linear spans, algebraic equations are handled via ideals generated by the polynomials they contain. We define such an ideal in the polynomial ring $\mathbb{R}[\bm{z}]$ as follows:
\begin{align*}
  \langle f_1, \ldots, f_r \rangle = \left\{\sum_{i=1}^r q_i f_i : q_i \in \mathbb{R}[\bm{z}]\right\} \subset \mathbb{R}[\bm{z}],
\end{align*}
where $f_1, \ldots, f_r \in \mathbb{R}[\bm{z}]$. The set $\langle f_1, \ldots, f_r \rangle$ is referred to as the ideal generated by $f_1,\ldots, f_r$.  The ideal $\langle f_1, \ldots, f_r \rangle$ is analogous to the span of a finite set of vectors. In each case, elements are formed through linear combinations, utilizing field coefficients for a span and polynomial coefficients for an ideal. We now show that $\langle f_1, \ldots, f_r \rangle$ satisfies the properties of an ideal in $\mathbb{R}[\bm{z}]$. 

\begin{lem}
  The set $\langle f_1, \ldots, f_r \rangle \subseteq \mathbb{R}[\bm{z}]$ is an ideal.
  \begin{proof}
    Let $\mathcal{J} = \langle f_1, \ldots, f_r \rangle$. Substituting the zero polynomials $h_1 = \cdots = h_r = 0$, yields $0 = \sum_{i=1}^r h_i f_i$, which implies $0 \in \mathcal{J}$.

    Consider $h, g \in \mathcal{J}$. By definition of $\mathcal{J}$, there exists $h = \sum_{i=1}^r h_i f_i, g = \sum_{i=1}^r g_i f_i$ such that $h_1,\ldots,h_r, g_1,\ldots,g_r \in \mathbb{R}[\bm{z}]$. Since 
    \begin{align*}
      h + g
      =
      \sum_{i=1}^r h_i f_i + \sum_{i=1}^r g_i f_i
      =
      \sum_{i=1}^r (h_i + g_i) f_i
    \end{align*}
    holds, we obtain $h + g \in \mathcal{J}$ by $h_1+g_1,\ldots,h_r+g_r \in \mathbb{R}[\bm{z}]$ and the definition of $\mathcal{J}$.

    Let $h \in \mathcal{J}$ and $c \in \mathbb{R}[\bm{z}]$. Subsequently, $h = \sum_{i=1}^r h_i f_i$ holds for some $h_1,\ldots,h_r \in \mathbb{R}[\bm{z}]$. Since
    \begin{align*}
      c h
      =
      c \sum_{i=1}^r h_i f_i
      =
      \sum_{i=1}^r (c h_i) f_i
    \end{align*}
    holds, we have $c h \in \mathcal{J}$ by $c h_1,\ldots,c h_r \in \mathbb{R}[\bm{z}]$ and the definition of $\mathcal{J}$. Thus, we conclude that $\mathcal{J} \subseteq \mathbb{R}[\bm{z}]$ is an ideal.
  \end{proof}
\end{lem}

In particular, an ideal generated by monomials is referred to as a monomial ideal, a concept utilized in the subsequent section. The following proposition establishes a property of monominial ideals \cite[Section 2.4, Lemmas 2 and 3]{CLO}. This property states that monomial ideal membership problems can be resolved without requiring specialized generator sets such as Gr\"obner bases; instead, membership can be determined by verifying whether the generators divide the monomials composing a given polynomial.

\begin{pro}\label{clo:2.4.2,2.4.3}
  Let $f = \sum_{\alpha = (\alpha_1,\ldots, \alpha_m) \in \mathbb{Z}_{\geq 0}^m} a_{\alpha} \bm{z}^{\alpha}$ be a polynomial in $\mathbb{R}[\bm{z}]$, and let $I = \langle \bm{z}^\beta : \beta \in B \rangle$ be a monomial ideal for some $B \subset \mathbb{Z}_{\geq 0}^m$. Then, $f \in I$ if and only if
  \begin{align*}
    \forall \alpha \in \mathbb{Z}_{\geq 0}^m \; \exists \beta \in B \; (a_{\alpha} \neq 0 \Rightarrow \bm{z}^{\beta} \mid \bm{z}^{\alpha})
  \end{align*}
\end{pro}

We introduce the notion of radicals, whose properties will be utilized in the subsequent section.

\begin{dfn}\label{def:rad}
  An ideal $I$ is radical if $f^{m} \in I$ for some integer $m \geq 1$ implies $f \in I.$
\end{dfn}

This section concludes with a fundamental result concerning radical ideals, known as Hilbert's Nullstellensatz \citep[Theorem 7.40]{BW}. This theorem establishes that the vanishing of a polynomial over a system of algebraic equations can be determined via radical ideal membership.

\begin{pro}\label{bw:7.40}
  Let $f$ be a polynomial in $\mathbb{R}[\bm{z}]$, and let $I = \langle p_1, \ldots, p_{\ell} \rangle$ be a radical ideal in $\mathbb{R}[\bm{z}]$. Then, $f \in I$ if and only if
  \begin{align*}
    \forall \bm{z} \in \mathbb{C}^m
    \left(
    \begin{cases}
      p_1(\bm{z}) = 0
      \\
      ~~ \vdots
      \\
      p_{\ell}(\bm{z}) = 0
    \end{cases}
    \Rightarrow f(\bm{z}) = 0
    \right)
  \end{align*}
\end{pro}

In the next section, we establish Theorem \ref{thm:thurstone-simple} by leveraging the properties of monomial ideals and radical ideals.

\section{Proof for Theorem \ref{thm:thurstone-simple}}\label{AppSec:TSS}

To establish Theorem \ref{thm:thurstone-simple}, it suffices to demonstrate that an initial solution satisfying Thurstone simple structure is generally not a stationary point of the orthomax criterion $Q_{\omega}(\Lambda) = Q_{\omega}(AT)$, where $0 \leq \omega \leq p$. The proof begins with the computation of the partial derivatives of the orthomax criteria. For $\lambda_{ij} = \sum_{l=1}^m a_{il} t_{lj}$, we derive the followings:
\begin{align*}
  \frac{\partial \lambda_{iv}^n}{\partial t_{uv}}
  =
  \frac{\partial}{\partial t_{uv}}
  \left[
    \sum_{l=1}^m a_{il} t_{lv}
  \right]^n
  =
  n a_{iu}
  \left[
    \sum_{l=1}^m a_{il} t_{lv}
  \right]^{n-1}
  = 
  n a_{iu} \lambda_{iv}^{n-1}
  ~~
  \mbox{and}
  ~~
  \frac{\partial \lambda_{ij}^n}{\partial t_{uv}} = 0. \quad (j \neq v)
\end{align*}
Consequently, we obtain
\begin{align}
  \frac{\partial Q_{\omega}}{\partial t_{uv}}
  &=
  \frac{\partial}{\partial t_{uv}}
  \left\{
  \sum_{i=1}^p \sum_{j=1}^m \lambda_{ij}^4 - \frac{\omega}{p}\sum_{j=1}^m \left(\sum_{i=1}^p \lambda_{ij}^2 \right)^2
  \right\}\notag
  \\
  &=
  \frac{\partial}{\partial t_{uv}}
  \left\{
  \sum_{i=1}^p \lambda_{iv}^4 - \frac{\omega}{p} \left(\sum_{i=1}^p \lambda_{iv}^2 \right)^2
  \right\}\notag
  \\
  &=
  \sum_{i=1}^p \frac{\partial \lambda_{iv}^4}{\partial t_{uv}} -
  \frac{\partial}{\partial t_{uv}}
  \left\{
  \frac{\omega}{p} \left(\sum_{i=1}^p \lambda_{iv}^2 \right)^2
  \right\}\notag
  \\
  &=
  \sum_{i=1}^p 4 a_{iu} \lambda_{iv}^3 - 2 \frac{\omega}{p} \left(\sum_{i=1}^p \lambda_{iv}^2 \right) \left(\sum_{i=1}^p \frac{\partial \lambda_{iv}^2}{\partial t_{uv}} \right)\notag
  \\
  &=
  \sum_{i=1}^p 4 a_{iu} \lambda_{iv}^3 - 2 \frac{\omega}{p} \left\| \lambda_{v} \right\|^2 \left(\sum_{i=1}^p 2 a_{iu} \lambda_{iv} \right)\notag
  \\
  &=
  4
  \sum_{i=1}^p
  a_{iu} \lambda_{iv} \left(\lambda_{iv}^2 - \frac{\omega}{p} \left\| \lambda_{v} \right\|^2 \right),
  \label{eq:tst}
\end{align}
where $\left\| \lambda_{v} \right\|^2 = \sum_{i=1}^p \lambda_{iv}^2$.

If an initial solution satisfying Thurstone simple structure is a stationary point of the orthomax criterion, then the substitution of the identity matrix into Eq. \eqref{eq:tst} must yield zero. Hence, we examine the following polynomial in the polynomial ring $\mathbb{R}[\bm{a}] = \mathbb{R}[a_{ij} : 1 \leq i \leq p, 1 \leq j \leq k]$:

\begin{align*}
  \frac{\partial Q_{\omega}}{\partial t_{uv}}(AI)
  &=
  4
  \left\{
  \sum_{i=1}^p
  a_{iu} a_{iv}
  \left(
  a_{iv}^2
  -
  \omega
  \left\| \lambda_{v} \right\|^2
  \right)
  \right\}.
\end{align*}

Given Proposition \ref{thm:eq7}, an initial solution is a stationary point of the orthomax criterion if and only if the following polynomial $f(\bm{a})$ vanishes for each $u \neq v$:
\begin{align*}
  f(\bm{a}) = 
  \sum_{i=1}^p
  a_{iu} a_{iv} \left(\left(a_{iu}^2 - a_{iv}^2 \right)- \omega\left(\left\| a_{u} \right\|^2-\left\| a_{v} \right\|^2\right)\right).
\end{align*}

By observing that the polynomial $f(\bm{a})$ depends solely on the elements of the $u$ and $v$ columns of an initial solution $A$, it follows that the vanishing of $f(\bm{a})$ is determined solely by the $u$ and $v$ columns of $A$. We assume that the initial solution satisfies a Thurstone simple structure of the class $(\gamma, \delta)$, where $\gamma + \delta \neq p$. Furthermore, we algebraically characterize the $u,v$ columns of the initial solution $A$ that conform to the Thurstone simple structure as follows:
\begin{enumerate}
\setcounter{enumi}{1}
\item\label{enu:Thurstone2:app}
  for each $j = u, v$
  \[
  \prod_{a \in \alpha} a = 0 \quad
  \text{for all } \alpha \subset A_j = \{ a_{1j}, \ldots, a_{pj} \} \text{ such that } |\alpha| = p - m + 1
  \]
\item\label{enu:Thurstone3:app}
  for each $i = 1, \ldots, \gamma_u$
  we have $a_{iu} = 0$,
  and
  for each $i = \gamma_u + 1, \ldots, \gamma_u + \gamma_v$
  we have 
  $a_{iv} = 0$
\item\label{enu:Thurstone4:app}
  for each $i = \gamma_u + \gamma_v + 1, \ldots, \gamma_u + \gamma_v + \delta$, we have $a_{iu} = a_{iv} = 0$.
\end{enumerate}
Note that conditions 3 and 4 do not lose generality through row swaps within column groups. Moreover, consideration of condition 1 of Thurstone's rule is unnecessary due to the following observations:
\begin{itemize}
\item
  rows satisfying condition 1 in the $u,v$ columns also fulfill conditions 3 or 4.
\item
  rows satisfying condition 1 in columns other than $u,v$ are unconstrained in the $u,v$ columns.
\end{itemize}
In other words, the initial solution $A$ satisfies Thurstone simple structure of class $(\gamma, \delta)$ if and only if the following algebraic equations are satisfied for each $u, v = 1, \ldots, k$, where $u \neq v$:
\begin{align}
  \begin{cases}
    0 = \prod_{a \in \alpha} a & \left(\alpha \subset A_u \mbox{ s.t. } |\alpha| = p - m + 1 \right)
    \\
    0 = \prod_{a \in \alpha} a & \left(\alpha \subset A_v \mbox{ s.t. } |\alpha| = p - m + 1 \right)
    \\
    0 = a_{iu} & \left(i = 1,\ldots,\gamma_u \right)
    \\
    0 = a_{iv} & \left(i = \gamma_u+1,\ldots,\gamma_u+\gamma_v \right)
    \\
    0 = a_{iu} = a_{iv} & \left(i = \gamma_u+\gamma_v+1,\ldots,\gamma_u+\gamma_v+\delta \right)
  \end{cases}.
  \label{initial-eq}
\end{align}

We want to show that if these equations are satisfied, the initial solution $A$ also is a sufficient condition for $A$ to be a stationary point. That is, we want to show that 
\begin{align}
  \bm{a} \in \mathbb{R}^{pk}
  \left(
  \eqref{initial-eq} \Rightarrow f(\bm{a}) = 0
  \right).
  \label{SHOW}
\end{align}
However, in the rest of this paper, we show that these conditions are generally not satisfied. Define
\begin{align*}
  F_u &=
  \left\{ \prod_{a \in \alpha} a : \alpha \subset A_u \mbox{ s.t. } |\alpha| = p - m + 1 \right\},
  \\
  F_v &=
  \left\{ \prod_{a \in \alpha} a : \alpha \subset A_v \mbox{ s.t. } |\alpha| = p - m + 1 \right\},
  \\
  H_u &=
  \left\{ a_{iu} : i = 1,\ldots,\gamma_u \right\},
  \\
  H_v &=
  \left\{ a_{iv} : i = \gamma_u+1,\ldots,\gamma_u+\gamma_v \right\},
  \\
  H_{uv} &=
  \left\{ a_{iu}, a_{iv} : i = \gamma_u+\gamma_v+1,\ldots,\gamma_u+\gamma_v+\delta \right\}.
\end{align*}
The ideal $\langle F_{u} \cup F_{v} \cup H_u \cup H_v \cup H_{uv} \rangle$ is a radical and monomial ideal in $\mathbb{R}[\bm{a}]$. While it is trivial that $\langle F_{u} \cup F_{v} \cup H_u \cup H_v \cup H_{uv} \rangle$ is a monomial ideal based on the definition, the following lemma demonstrates that it is also a radical ideal. 

\begin{lem}
The ideal $\langle F_{u} \cup F_{v} \cup H_u \cup H_v \cup H_{uv} \rangle$ is a radical ideal in $\mathbb{R}[\bm{a}]$.
\begin{proof}
Assume for contradiction that $\langle F_{u} \cup F_{v} \cup H_u \cup H_v \cup H_{uv} \rangle$ is not a radical ideal. Then, there exists $f \in \mathbb{R}[\bm{a}]$ such that $f^m \in \langle F_{u} \cup F_{v} \cup H_u \cup H_v \cup H_{uv} \rangle$ and $f \not\in \langle F_{u} \cup F_{v} \cup H_u \cup H_v \cup H_{uv} \rangle$ for some $m \geq 1$, by Definition \ref{def:rad}. Let $f = \sum_{\alpha \in \mathbb{Z}_{\geq 0}^{pk}} c_{\alpha} \bm{a}^{\alpha}, c_{\alpha} \in \mathbb{R}$. 

Since $\langle F_{u} \cup F_{v} \cup H_u \cup H_v \cup H_{uv} \rangle$ is a monomial ideal, we invoke Proposition \ref{clo:2.4.2,2.4.3}, which provides the following equivalence condition:
\begin{align*}
&
\begin{cases}
f^m \in \langle F_{u} \cup F_{v} \cup H_u \cup H_v \cup H_{uv} \rangle
\\
f \not\in \langle F_{u} \cup F_{v} \cup H_u \cup H_v \cup H_{uv} \rangle
\end{cases}
\\
\Longleftrightarrow
&
\begin{cases}
\forall t_1 \in \mbox{Mon}(f^m)\, \exists s_1 \in F_{u} \cup F_{v} \cup H_u \cup H_v \cup H_{uv} \; (s_1 \mid t_1)
\\
\forall t_2 \in \mbox{Mon}(f) \, \forall s_2 \in F_{u} \cup F_{v} \cup H_u \cup H_v \cup H_{uv} \; (s_2 \nmid t_2)
\end{cases},
\tag{$\ast$}
\end{align*}
where $\mbox{Mon}(f^m)$ and $\mbox{Mon}(f)$ denote the sets of monomials of $f^m$ and $f$, respectively.

It follows that any monomial $t_1 \in \mbox{Mon}(f^m)$ is a product of elements from $\mbox{Mon}(f)$, thereby contradicting $(\ast)$.
\end{proof}
\end{lem}

Moreover, some monomials contained in the polynomial $r$ which is a component of $f$ are not divisible by any monomials in $F_{u} \cup F_{v} \cup H_u \cup H_v \cup H_{uv}$:
\begin{align*}
  r(\bm{a}) =
  \sum_{i=k+\ell+1}^p
  a_{iu} a_{iv}
  \left(
  \left(
  a_{iu}^2
  -
  a_{iv}^2
  \right)
  -
  \omega
  \left(
  \left\| a_{u} \right\|^2-\left\| a_{v} \right\|^2
  \right)
  \right)
\end{align*}
Hence, as the relation \eqref{SHOW} is generally not satisfied by Propositions \ref{clo:2.4.2,2.4.3} and \ref{bw:7.40}, Theorem \ref{thm:thurstone-simple} is established.

% \section{Details for Theorem \ref{thm:BH}}\label{AppSec:BH}

\end{document}